\documentclass[aos,MSNbibl,nameyear,dvips]{arximspdf}
\usepackage{graphicx}
\doi{10.1214/11-AOS943}
\volume{39}
\issue{6}
\pubyear{2011}
\firstpage{3290}
\lastpage{3319}

\makeatletter

\def\bsuffix #1{#1}

\newcommand{\dddot}{\hspace*{-0.5pt}\dot{\hphantom{z}}\hspace*{-1.5pt}\ddot{\hphantom{0}}\hspace*{-7pt}z}

\newcommand{\eps}{\varepsilon}
\newcommand{\R}{\mathbb{R}}
\newcommand{\mH}{\mathcal{H}}

\renewcommand{\P}{\mathrm{P}}
\newcommand{\E}{\mathrm{E}}
\newcommand{\var}{\operatorname{var}}
\newcommand{\FWER}{\operatorname{FWER}}
\newcommand{\FDR}{\operatorname{FDR}}
\newcommand{\Power}{\operatorname{Power}}
\newcommand{\Bon}{\operatorname{Bon}}
\newcommand{\BH}{\operatorname{BH}}
\newcommand{\argmax}{\mathop{\arg\max}}
\newcommand{\SNR}{\operatorname{SNR}}

\newtheorem{theorem}{Theorem}
\newtheorem{lemma}[theorem]{Lemma}
\newtheorem{prop}[theorem]{Proposition}

\newproclaim{alg}[theorem]{Algorithm}
\newproclaim{example}[theorem]{Example}

\makeatother

\begin{document}
\begin{frontmatter}

\title{Multiple testing of local maxima for detection of~peaks in 1D}
\runtitle{Multiple testing of local maxima in 1D}

\begin{aug}
\author[A]{\fnms{Armin} \snm{Schwartzman}\corref{}\thanksref{t1}\ead[label=e1]{armins@hsph.harvard.edu}},
\author[A]{\fnms{Yulia} \snm{Gavrilov}\thanksref{t1}\ead[label=e2]{yuliagavrilov@gmail.com}}
\and
\author[B]{\fnms{Robert J.} \snm{Adler}\thanksref{t2}\ead[label=e3]{robert@ee.technion.ac.il}}
\runauthor{A. Schwartzman, Y. Gavrilov and R. J. Adler}
\affiliation{Harvard School of Public Health, Harvard School of Public Health
and~Technion, Israel Institute of Technology}
\address[A]{A. Schwartzman\\
Y. Gavrilov\\
Department of Biostatistics \\
Harvard School of Public Health \\
and \\
Dana-Farber Cancer Institute \\
450 Brookline Ave., CLS 11007 \\
Boston, Massachusetts 02446 \\
USA\\
\printead{e1} \\
\hphantom{E-mail: }\printead*{e2}}
\address[B]{R. J. Adler\\
Department of Electrical Engineering \\
Technion, Israel Institute of Technology \\
Haifa 32000\\
Israel \\
\printead{e3}} 
\end{aug}

\thankstext{t1}{Supported in part by the Claudia Adams Barr Program in
Cancer Research,
the William F. Milton Fund and NIH Grant P01-CA134294.}

\thankstext{t2}{Supported in part by US--Israel Binational Science
Foundation, 2008262 and Israel Science Foundation, 853/10.}

\received{\smonth{5} \syear{2011}}
\revised{\smonth{11} \syear{2011}}

%
\begin{abstract}
A topological multiple testing scheme for one-dimensional domains is
proposed where, rather than testing every spatial or temporal location
for the presence of a signal, tests are performed only at the local
maxima of the smoothed observed sequence. Assuming unimodal true peaks
with finite support and Gaussian stationary ergodic noise, it is shown
that the algorithm with Bonferroni or Benjamini--Hochberg correction
provides asymptotic strong control of the family wise error rate and
false discovery rate, and is power consistent, as the search space and
the signal strength get large, where the search space may grow
exponentially faster than the signal strength. Simulations show that
error levels are maintained for nonasymptotic conditions, and that
power is maximized when the smoothing kernel is close in shape and
bandwidth to the signal peaks, akin to the matched filter theorem in
signal processing. The methods are illustrated in an analysis of
electrical recordings of neuronal cell activity.
\end{abstract}

%

\end{frontmatter}

\section{Introduction}

One of the most challenging aspects of multiple testing problems in
spatial and temporal domains is how to account for the spatial or
temporal structure in the underlying signal. The usual paradigm
considers a separate test at each observed location. However, the
interest is usually in detecting signal regions that span several
neighboring locations. This paper considers a new multiple testing
paradigm for spatial and temporal domains where tests are not performed
at every observed location, but only at the local maxima of the
observed data, seen as representatives of underlying signal peak
regions. The proposed inference is not pointwise but topological, based
on the observed local maxima as topological features.

In pointwise testing, the control of family-wise error rate (FWER), now
common in neuroimaging, was established by Keith Worsley
[\citet{Taylor2007}, Worsley et~al.
(\citeyear{Worsley1996a,Worsley2004})], who exploited the Euler
characteristic heuristic for approximating the distribution of the
maximum of a random field [\citet{Adler2007},
\citet{Adler2010}]. Methods for controlling the false discovery
rate (FDR) [\citet{Benjamini1995}] are also applied routinely in
this setting, but the spatial structure is difficult to incorporate and
often ignored [\citet{Genovese2002}, \citet{Nichols2003},
\citet{Schwartzman2008a}].

Despite pointwise testing being so common, the real interest is usually
not in detecting individual locations, but connected regions or
clusters. This has prompted the adaptation of discrete FDR methods to
pre-defined clusters [\citet{Heller2007},
\citet{Heller2006}], and the use of Gaussian random field theory
for computing $p$-values corresponding to the height, extent and mass of
clusters obtained by pre-thresholding the observed random field
[\citet{Poline1997}, \citet{Zhang2009}]. Perone~Pacifico
et~al. (\citeyear{Pacifico2004,Pacifico2007}) proposed data-dependent
thresholds so that FDR is controlled at the cluster level, using
Gaussian random field theory to approximate the null distribution.
However, the definition of Type I error for clusters requires a
tolerance parameter for the overlap between a discovered cluster and
the null region [\citet{Pacifico2004}], while spatial smoothing,
which is often applied for improving signal-to-noise ratio (SNR),
creates the need to remove the spread of the signal over the null
region to avoid error inflation [\citet{Pacifico2007}].
\citet{Chumbley2009} have argued that current cluster methods are
unsatisfactory because, just like marginal FDR procedures, they rely on
the basic premise of having a test at each spatial location; instead,
inference should be topological.

This article proposes a different multiple testing paradigm where tests
are performed, not at each spatial or temporal location, but only at
the local maxima of the smoothed data, seen as topological
representatives of their neighborhood region or cluster. A similar idea
was recently proposed independently by \citet{Chumbley2010}, but they
did not consider whether Type I error could be controlled. Here we
extend the classical control of FWER via the global maximum to control
of both FWER and FDR via local maxima. Because the distributional
theory for local maxima of random fields is more difficult than that
for global maxima, this paper only considers one-dimensional domains
(spatial or temporal), where closed-form solutions exist, leaving the
two- and three-dimensional cases for future work.

\begin{figure}

\includegraphics{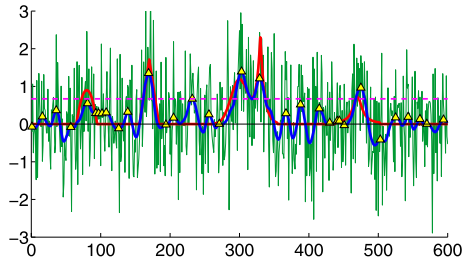}

\caption{Simulated observed sequence $y(t)$ (green) and
smoothed sequence $y_\gamma(t)$ (blue) over five underlying true peaks
of different shapes comprising $\mu(t)$ (red). Out of 33 local maxima
of $y_\gamma(t)$ (yellow), the BH detection threshold at FDR level 0.2
(dashed magenta) selects five, one of which is a false positive. At
this noise level, four out of five true peaks are detected. Note that
this bandwidth is able to distinguish the overlapping peaks.}
\label{figsimulexample}
\end{figure}

Our general proposed algorithm consists of the following steps:
\begin{longlist}
\item[(1)] \textit{Kernel smoothing}: to increase SNR
[\citet{Smith2009}, \citet{Worsley1996b}].
\item[(2)] \textit{Candidate peaks}: find local maxima of the smoothed sequence.
\item[(3)] \textit{$p$-values}: computed at each local maximum under the
complete null hypothesis of no signal anywhere.
\item[(4)] \textit{Multiple testing}: apply a multiple testing procedure
and declare as detected peaks those local maxima whose $p$-values are significant.
\end{longlist}
In this paper, the $p$-values in step (3) are computed using theory of
Gaussian processes. For step (4), we consider two standard multiple
testing procedures: Bonferroni to control FWER and Benjamini--Hochberg
(BH) [\citet{Benjamini1995}] to control FDR. The algorithm is
illustrated by a simulated example in Figure \ref{figsimulexample}.

%

We study the theoretical properties of the above algorithm under a
specific signal-plus-noise model and then relax these assumptions in
the simulations. For Type I errors to be well defined, the signal is
modeled as if composed of unimodal peak regions, each considered
detected if a significant local maximum occurs inside its finite
support. For simplicity, we concentrate on positive signals and
one-sided tests, but this is not crucial. For tractability, the theory
assumes that the observation noise follows a smooth stationary ergodic
Gaussian process. This assumption permits an explicit formula for
computing the $p$-values corresponding to local maxima of the observed
process. The distribution of the height of a local maximum of a
Gaussian process is not Gaussian but has a heavier tail, and its
computation requires careful conditioning based on the calculus of Palm
probabilities [\citet{Adler2010}, \citet{Cramer1967}].

An interesting and challenging aspect of inference for local maxima is
the fact that the number of tests, equal to the number of observed
local maxima, is random. The multiple testing literature usually
assumes the number of tests to be fixed. We overcome this difficulty
with an asymptotic argument for large search space, so that by
ergodicity, the error behaves approximately as it would if the number
of tests were equal to its expected value.


In order to achieve strong control of FWER and FDR, the asymptotics for
large search space are combined with asymptotics for strong signal. The
strong signal assumption asymptotically eliminates the false positives
caused by the smoothed signal spreading into the null regions, by
assuring that each signal peak region is represented by only one
observed local maximum within the true domain with probability tending
to one. The strong signal assumption is not restrictive in the sense
that the search space may grow exponentially faster. Simulations show
that error levels are maintained at finite search spaces and moderate
signal strength.

Defining detection power as the expected fraction of true peaks
detected, we prove that the algorithm is consistent in the sense that
its power tends to one under the above asymptotic conditions. We find
that the optimal smoothing kernel is approximately that which is
closest in shape and bandwidth to the signal peaks to be detected, akin
to the so-called matched filter theorem in signal processing
[\citet{Pratt1991}, \citet{Simon1995}]. This optimal bandwidth is much larger
than the usual optimal bandwidth for nonparametric regression.

In one dimension, the problem of identifying significant local maxima
is similar to that of peak detection in signal processing [e.g.,
\citet{Arzeno2008}, \citet{Baccus2002}, \citet{Brutti2005},
\citet{Harezlak2008}, \citet{Morris2005}, \citet{Yasui2003}].
In this literature, though large, the detection threshold is
predominantly chosen heuristically and conservatively. Our multiple
testing viewpoint provides a formal mechanism for choosing the
detection threshold, allowing detection under higher noise conditions.
This view also eliminates the need to estimate an unknown number of
peak location parameters, encountered in the signal estimation approach
[Li and Speed (\citeyear{Li2000,Li2004}), \citet{OBrien1994},
\citet{Tibshirani2005}].

We illustrate our procedure with a data set of neural electrical
recordings, where the objective is to detect action potentials
representing cell activity [\citet{Baccus2002},
Segev et al. (\citeyear{Segev2004})]. The noise parameters and signal peak shape are
estimated from a training set and then applied to a~test set for peak
detection.


The data analysis and all simulations were implemented in \texttt{Matlab}.

\section{Theory}
\label{sectheory}
\subsection{The model}
\label{secmodel}
Consider the signal-plus-noise model
%
\begin{equation}
\label{eqsignal+noise}
y(t) = \mu(t) + z(t),\qquad t \in\R,
\end{equation}
where the signal $\mu(t)$ is a train of unimodal positive peaks of the form
%
\begin{equation}
\label{eqmu}
\mu(t) = \sum_{j=-\infty}^\infty a_j h_j(t),\qquad a_j > 0,
\end{equation}
and the peak shape $h_j(t) \ge0$ has compact connected support $S_j =
\{t\dvtx h_j(t) > 0\}$ and unit action $\int_{S_j} h_j(t) \,dt = 1$ for
each $j$. Let $w_\gamma(t) \ge0$ with bandwidth parameter $\gamma>0$
be a unimodal kernel with compact connected support and unit action.
Convolving the process (\ref{eqsignal+noise}) with the kernel
$w_\gamma(t)$ results in the smoothed process
%
\begin{equation}
\label{eqconv}
y_\gamma(t) = w_\gamma(t) * y(t) =
\int_{-\infty}^\infty w_\gamma(t-s) y(s) \,ds = \mu_\gamma(t) +
z_\gamma(t),
\end{equation}
where the smoothed signal and smoothed noise are defined as
%
\begin{equation}
\label{eqmu-gamma}
\mu_\gamma(t) = w_\gamma(t) * \mu(t) = \sum_{j=-\infty}^\infty
a_j h_{j,\gamma}(t),\qquad z_\gamma(t) = w_\gamma(t) * z(t).
\end{equation}

For each $j$, the smoothed peak shape $h_{j,\gamma}(t) = w_\gamma
(t)*h_j(t) \ge0$ is unimodal and has compact connected support
$S_{j,\gamma}$ and unit action. For each $j$, we require that
$h_{j,\gamma}(t)$ is twice differentiable in the interior of
$S_{j,\gamma}$ and has no other critical points within its support.
For simplicity, the theory requires that the supports $S_{j,\gamma}$
do not overlap (but this is not required in practice, as shown via
simulations in Section \ref{secsimulations}). The smoothed noise
$z_\gamma(t)$ defined by (\ref{eqconv}) and (\ref{eqmu-gamma}) is
assumed to be a zero-mean thrice differentiable stationary ergodic
Gaussian process.

\subsection{The STEM algorithm}
\label{secalg}
Suppose we observe $y(t)$ defined by (\ref{eqsignal+noise}) in the
segment $[-L/2,L/2]$, which contains $J$ peaks. We call the following
procedure STEM (Smoothing and TEsting of Maxima).
%
\begin{alg}[(STEM algorithm)]
\label{algSTEM}

(1) \textit{Kernel smoothing}:
construct the process (\ref{eqconv}), ignoring the boundary effects at
$\pm L/2$.\vspace*{-5pt}
\begin{longlist}[(2)]
\item[(2)] \textit{Candidate peaks}:
find the set of local maxima of $y_\gamma(t)$ in $[-L/2,L/2]$
%
\begin{equation}
\label{eqT}
\tilde{T} = \biggl\{ t \in\biggl[-\frac{L}{2},\frac{L}{2}\biggr]\dvtx
\dot{y}_{\gamma}(t) = \frac{dy_{\gamma}(t)}{dt} = 0,
\ddot{y}_{\gamma}(t) = \frac{d^2 y_{\gamma}(t)}{dt^2} < 0 \biggr\}.
\end{equation}
\item[(3)] \textit{$p$-values}:
for each $t \in\tilde{T}$ compute the $p$-value $p_\gamma(t)$ for
testing the (conditional) hypothesis
\[
\mH_{0}(t)\dvtx \mu(t) = 0 \quad\mbox{vs.}\quad \mH_{A}(t)\dvtx \mu
(t) > 0,\qquad
t \in\tilde{T}.
\]
\item[(4)] \textit{Multiple testing}:
let $\tilde{m}$ be the number of tested hypotheses, equal to the
number of local maxima in $\tilde{T}$. Apply a multiple testing
procedure on the set of $\tilde{m}$ $p$-values $\{p_\gamma(t), t \in
\tilde{T}\}$, and declare significant all peaks whose $p$-values are
smaller than the significance threshold.\vadjust{\goodbreak}
\end{longlist}
\end{alg}

Steps (1) and (2) above are well defined under the model assumptions (for
data on a grid, local maxima are defined as points higher than their
neighbors). Step (3) is detailed in Section \ref{secpvalue} below. For
step (4), we use the Bonferroni procedure to control FWER and the BH
procedure to control FDR. To apply Bonferroni at level~$\alpha$,
declare significant all peaks whose $p$-values are smaller than $\alpha
/\tilde{m}$. To apply BH at level $\alpha$, find the largest index~$k$
for which the $i$th smallest $p$-value is smaller than $i\alpha
/\tilde{m}$, and declare as significant the $k$ peaks with smallest
$p$-values. Notice that, in contrast to the usual application of the
Bonferroni and BH procedures, the number of tests $\tilde{m}$ is random.

\subsection{$p$-values}
\label{secpvalue}
Given the observed heights $y_\gamma(t)$ at the local maxima \mbox{$t\in
\tilde{T}$}, the $p$-values in step (3) of Algorithm \ref{algSTEM} are
computed as
%
\begin{equation}
\label{eqp-value}
p_\gamma(t) = F_\gamma[y_\gamma(t)],\qquad t\in\tilde{T},
\end{equation}
where
%
\begin{equation}
\label{eqpalm}
F_\gamma(u) = \P\{z_\gamma(t) > u | t \in\tilde{T}
\}
\end{equation}
denotes the \textit{right} cumulative distribution function (cdf) of
$z_\gamma(t)$ at the local maxima $t \in\tilde{T}$, evaluated under
the complete null hypothesis $\mu(t) = 0, \forall t$.

The conditional distribution (\ref{eqpalm}) is called a Palm
distribution [\citet{Adler2010}, Chapter 6]. Unlike the marginal
distribution of $z_\gamma(t)$, it is not Gaussian but stochastically
greater. This is because the point of evaluation $t\in\tilde{T}$ is not
a fixed point $t\in\R$, but the random location of a local maximum of
$z_\gamma(t)$. Moreover, the conditioning event has probability zero.
The Palm distribution (\ref{eqpalm}) has a closed-form expression,
originally obtained by Cram{\'e}r and
Leadbetter [(\citeyear{Cramer1967}), Chapter 11] (equation
11.6.14), using the well-known Kac--Rice formula [\citet{Rice1945},
\citet{Adler2007}, Chapter 11]. A direct application, borrowing
notation from those sources, gives the following.
%
\begin{prop}
\label{propPalm}
Suppose the assumptions of Section \ref{secmodel} hold and that $\mu
(t) = 0, \forall t$. Define the moments
%
\begin{equation}
\label{eqmoments}
\sigma^2_\gamma= \var[z_\gamma(t)],\qquad
\lambda_{2,\gamma} = \var[\dot{z}_\gamma(t)],\qquad
\lambda_{4,\gamma} = \var[\ddot{z}_\gamma(t)].
\end{equation}
Then the distribution (\ref{eqpalm}) is given by
%
\begin{equation}
\label{eqdistr}
F_\gamma(u) = 1 - \Phi\Biggl(u \sqrt{\frac{\lambda_{4, \gamma
}}{\Delta}}\Biggr) + \sqrt{\frac{2\pi\lambda^2_{2, \gamma
}}{\lambda_{4, \gamma}\sigma^2_\gamma}}\phi\biggl(\frac{u}{\sigma
_\gamma}\biggr)\Phi\Biggl( u \sqrt{\frac{\lambda^2_{2, \gamma
}}{\Delta\sigma_\gamma^2}}\Biggr),
\end{equation}
where $\Delta= \sigma^2_\gamma\lambda_{4,\gamma} - \lambda
_{2,\gamma}^2$, and $\phi(x)$, $\Phi(x)$ are the standard normal
density and cdf, respectively.
\end{prop}

The quantities $\sigma^2_\gamma$, $\lambda_{2,\gamma}$ and $\lambda
_{4,\gamma}$ in Proposition \ref{propPalm} depend on the kernel~%
$w_\gamma(t)$ and the autocorrelation\vadjust{\goodbreak} function of the original noise
process $z(t)$. Explicit expressions may be obtained, for instance, for
the following Gaussian autocorrelation model, which we use later in the
simulations.
%
\begin{example}[(Gaussian autocorrelation model)]
\label{exGaussian-ACVF}
Let the noise $z(t)$ in (\ref{eqsignal+noise}) be constructed as
\[
z(t) = \sigma\int_{-\infty}^\infty\frac{1}{\nu} \phi\biggl(\frac
{t-s}{\nu}\biggr) \,dB(s),\qquad \sigma, \nu> 0,
\]
where $B(s)$ is standard Brownian motion and $\nu> 0$. Convolving with
a~Gaussian kernel $w_\gamma(t) = (1/\gamma)\phi(t/\gamma)$ with
$\gamma> 0$ as in (\ref{eqmu-gamma}) produces a zero-mean infinitely
differentiable stationary ergodic Gaussian process
\[
z_\gamma(t) = w_\gamma(t) * z(t) = \sigma\int_{-\infty}^\infty
\frac{1}{\xi} \phi\biggl(\frac{t-s}{\xi}\biggr) \,dB(s),\qquad
\xi= \sqrt{\gamma^2 + \nu^2},
\]
with moments (\ref{eqmoments}) given by
$\sigma^2_\gamma= \sigma^2/(2\sqrt{\pi} \xi)$,
$\lambda_{2,\gamma} = \sigma^2/(4\sqrt{\pi} \xi^3)$,
$\lambda_{4,\gamma} = 3\sigma^2/(8\sqrt{\pi} \xi^5)$.
The above expressions may be used as approximations if the kernel,
required to have finite support, is truncated at $t = \pm\gamma d$ for
moderately large~$d$, say $d=3$.
\end{example}

\subsection{Error definitions}
\label{secerrors}
Because truly detected peaks may be shifted with respect to the true
peaks as a result of noise, we define a significant local maximum to be
a true positive if it falls anywhere inside the support of a~true peak.
Conversely, we define it to be a false positive if it falls outside the
support of any true peak. Assuming the model of Section \ref
{secmodel}, define the signal region $\mathbb{S}_1$ and null region
$\mathbb{S}_0$, respectively, by
%
\begin{equation}
\label{eqnull-region}
\mathbb{S}_1 = \bigcup_{j=1}^J S_j
\quad\mbox{and}\quad
\mathbb{S}_0 = \biggl[-\frac{L}{2},\frac{L}{2}\biggr] \Bigm\backslash
\Biggl( \bigcup_{j=1}^J S_j \Biggr).
\end{equation}
For a significance threshold $u$, the total number of detected peaks
and the number of falsely detected peaks are
\[
R(u) = \#\{t\in\tilde{T} \dvtx y_\gamma(t) > u\} \quad\mbox{and}\quad
V(u) = \#\{t\in\tilde{T}\cap\mathbb{S}_0 \dvtx y_\gamma(t) >u\},
\]
respectively. Both are defined as zero if $\tilde{T}$ is empty. The
FWER is defined as the probability of obtaining at least one falsely
detected peak
%
\begin{equation}
\label{eqFWER}\qquad
\FWER(u) = \P\{ V(u)\ge1 \} = \P\Bigl\{ \tilde
{T}\cap\mathbb{S}_0 \ne\varnothing\mbox{ and }
\max_{t \in\tilde{T}\cap\mathbb{S}_0} y_\gamma(t) > u \Bigr\}.
\end{equation}
The FDR is defined as the expected proportion of falsely detected peaks
%
\begin{equation}
\label{eqFDR}
\FDR(u) = \E\biggl\{ \frac{V(u)}{R(u)\vee1} \biggr\}.
\end{equation}

\begin{figure}

\includegraphics{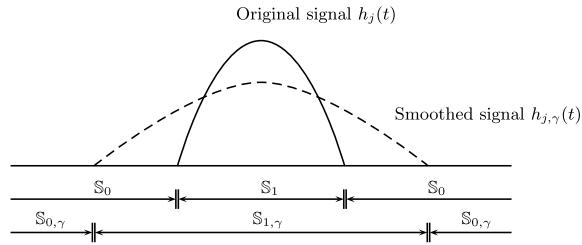}

\caption{Schematic signal and null regions, before and after smoothing,
in the vicinity of one signal peak.}
\label{figs0}
\end{figure}

Note that the above definitions are with respect to the original signal
support $\mathbb{S}_1$, while the inference is carried out using the
smoothed observed process $y_\gamma(t)$. Kernel smoothing enlarges the
signal support and increases the probability of obtaining false
positives in the null regions neighboring the signal
[\protect\citet{Pacifico2007}]. In contrast to (\protect\ref{eqnull-region}), the
smoothed signal region $\mathbb{S}_{1,\gamma} \supset\mathbb{S}_1$ and
smoothed null region $\mathbb{S}_{0,\gamma} \subset\mathbb{S}_0$ are
%
\begin{equation}
\label{eqnull-region-smoothed}
\mathbb{S}_{1, \gamma} = \bigcup_{j=1}^J S_{j, \gamma}
\quad\mbox{and}\quad
\mathbb{S}_{0, \gamma} = \biggl[-\frac{L}{2},\frac{L}{2}\biggr]
\Bigm\backslash\Biggl( \bigcup_{j=1}^J S_{j, \gamma}\Biggr),
\end{equation}
respectively (Figure \ref{figs0}). We call the difference between the
expanded signal support and the true signal support the transition region
%
\begin{equation}
\label{eqtransition-region}
\mathbb{T}_{\gamma} = \mathbb{S}_{1, \gamma}\setminus\mathbb{S}_1
=\mathbb{S}_{0}\setminus\mathbb{S}_{0, \gamma} = \bigcup_{j=1}^J
T_{j, \gamma},
\end{equation}
where $T_{j, \gamma} = S_{j, \gamma} \setminus S_j$ is the transition
region corresponding to each peak $j$.

In general, a true peak may produce more than one significant local
maximum, affecting the interpretation of definition (\ref{eqFDR}) and
the nonasymptotic validity of the FDR controlling procedure. However,
as explained below, this multiplicity is unlikely to occur for strong
signals, assuring validity at least asymptotically under that regime.
The simulations of Section \ref{secsim1} show it not to be problematic
in nonasymptotic situations for moderate signals and appropriate smoothing.

\subsection{Strong control of FWER}
\label{secFWER}
In Algorithm \ref{algSTEM}, step (3) produces a list of $\tilde{m}$
$p$-values. If the Bonferroni correction is applied in step (4) with
level\
$\alpha\in(0,1)$, then the null hypothesis $\mH_0(t)$ at $t \in
\tilde{T}$ is rejected if
%
\begin{equation}
\label{eqthresh-Bon-random}
p_\gamma(t) < \frac{\alpha}{\tilde{m}}
\quad\iff\quad y_\gamma(t) >
\tilde{u}_{\Bon} = F_\gamma^{-1} \biggl(\frac{\alpha}{\tilde
{m}}\biggr),
\end{equation}
where $\alpha/\tilde{m}$ is defined as 1 if $\tilde{m}=0$.
Recall that, in contrast to the usual Bonferroni algorithm, the number
of $p$-values $\tilde{m}$ is random.

Define the conditions:
\begin{longlist}[(C2)]
\item[(C1)] The assumptions of Section \ref{secmodel} hold.
\item[(C2)] $L \to\infty$ and $a = \inf_j a_j \to\infty$, such
that $(\log L)/a^2 \to0$ and \mbox{$J/L \to A_1$} with $0 < A_1 < 1$.
\end{longlist}
%
\begin{theorem}
\label{thmFWER}
Suppose that Algorithm \ref{algSTEM} is applied with the Bonferroni
threshold $\tilde{u}_{\Bon}$ (\ref{eqthresh-Bon-random}). Then,
under conditions \textup{(C1)} and \textup{(C2)},
\[
\limsup\FWER(\tilde{u}_{\Bon}) \le\alpha.
\]
\end{theorem}

The proof of Theorem \ref{thmFWER} is given in Section \ref{appFWER}.
The large search space assumption in (C2) solves the problem of $\tilde
{m}$ being random, implying that by the weak law of large numbers, the
ratio $\tilde{m}/L$ is close to its expectation $\E[\tilde{m}/L]$
for large $L$. Thus the Bonferroni procedure with random threshold~(\ref{eqthresh-Bon-random}) has asymptotically the same error control
properties as if the threshold were deterministic and equal to
%
\begin{equation}
\label{eqthresh-Bon-fixed}
u_{\Bon}^* = F_\gamma^{-1}\biggl(\frac{\alpha}{\E[\tilde
{m}]}\biggr) \approx F_\gamma^{-1}\biggl(\frac{\alpha/L}{A_1 + \E
[\tilde{m}_{0,\gamma}(0,1)]}\biggr),
\end{equation}
where
%
\begin{equation}
\label{eqexpected-local-maxima}
\E[\tilde{m}_{0,\gamma}(0,1)] = \frac{1}{2\pi}\sqrt{\frac
{\lambda_{4,\gamma}}{\lambda_{2,\gamma}}}
\end{equation}
is the expected number of local maxima of $z_\gamma(t)$ in the unit
interval $(0,1)$ [\citet{Cramer1967}, Chapter 10].

The strong signal assumption in (C2) implies (Lemma \ref
{lemmaunique-max} in Section \ref{applemmas}) that, with probability
tending to 1, no local maxima are obtained in the transition region
$\mathbb{T}_{\gamma}$ (\ref{eqtransition-region}), and exactly one
local maxima is obtained for each signal peak in $\mathbb{S}_1$. This
avoids the error inflation due to smoothing and provides the approximation
in~(\ref{eqthresh-Bon-fixed}). The proof of Lemma \ref
{lemmaunique-max} shows that the asymptotic rates are exponential and
controlled partially by the smallest absolute derivative of the
smoothed peak shape in the transition region and the curvature of the
smoothed peak shape at the mode.


\subsection{Control of FDR}
\label{secFDR-control}
Suppose the BH procedure is applied in step (4) of Algorithm \ref
{algSTEM}. For a fixed $\alpha\in(0,1)$, let $k$ be the largest index
for which the $i$th smallest $p$-value is less than $i\alpha/\tilde
{m}$. Then the null hypothesis $\mH_0(t)$ at $t \in\tilde{T}$ is
rejected if
%
\begin{equation}
\label{eqthresh-BH-random}
p_\gamma(t) < \frac{k\alpha}{\tilde{m}}
\quad\iff\quad y_\gamma(t) >
\tilde{u}_{\BH} = F_\gamma^{-1} \biggl(\frac{k\alpha}{\tilde
{m}}\biggr),
\end{equation}
where $k\alpha/\tilde{m}$ is defined as 1 if $\tilde{m}=0$.
%
\begin{theorem}
\label{thmFDR}
Suppose that Algorithm \ref{algSTEM} is applied with the BH threshold~$\tilde{u}_{\BH}$~(\ref{eqthresh-BH-random}). Then, under conditions
\textup{(C1)} and \textup{(C2)},
\[
\limsup\FDR(\tilde{u}_{\BH}) \le\alpha.
\]
\end{theorem}

The proof of Theorem \ref{thmFDR} is given in Section \ref{appFDR}.
The asymptotic assumptions~(C2), imply that the BH procedure with
random threshold (\ref{eqthresh-BH-random}) has asymptotically the
same error control properties as if the threshold were deterministic
and equal to
%
\begin{equation}
\label{eqthresh-BH-fixed}
u^*_{\BH} = F_\gamma^{-1}\biggl(\frac{\alpha A_1}{A_1 + \E[\tilde
{m}_{0,\gamma}(0,1)](1-\alpha)}\biggr),
\end{equation}
where $\E[\tilde{m}_{0,\gamma}(0,1)]$ is given by
(\ref{eqexpected-local-maxima}). The threshold
(\ref{eqthresh-BH-random}) can be viewed as the largest solution of
the equation $\alpha G(u) = F_\gamma(u)$, where $G(u)$ is the empirical
right cumulative distribution function of $y_\gamma(t), t\in
\tilde{T}$ [\citet{Genovese2002}]. Taking the limit of that equation as
$L$ gets large yields the solution (\ref{eqthresh-BH-fixed}).

As before, the strong signal assumption in (C2) implies that there
exists exactly one significant local maximum at each true peak with
probability tending to 1 (Lemma \ref{lemmaunique-max} in Section \ref
{applemmas}), avoiding error inflation in the transition region and
justifying the interpretation of definition (\ref{eqFDR}) as the
expected proportion of falsely discovered peaks. Again, the asymptotic
rates are exponential and controlled partially by the smallest absolute
derivative of the smoothed peak shape in the transition region and the
curvature of the smoothed peak shape at the mode.

Notice that, in contrast to the asymptotic Bonferroni threshold~$u^*_{\Bon}$~(\ref{eqthresh-Bon-fixed}) which grows unbounded with
increasing $L$, the asymptotic BH thresh-\break old~$u^*_{\BH}$~(\ref
{eqthresh-BH-fixed}) is finite.

\subsection{Power}
\label{secpower}
Recall from Section \ref{secerrors} that a significant local maximum
is considered a true positive if it falls in the true signal region
$\mathbb{S}_{1}$.
We define the power of Algorithm \ref{algSTEM} as the expected
fraction of true discovered peaks
%
\begin{eqnarray}
\label{eqpower}
\Power(u) &=& \E\Biggl[ \frac{1}{J} \sum_{j=1}^{J} 1\Bigl(
\tilde{T} \cap S_j \ne\varnothing\mbox{ and } \max_{\tilde{t}
\in\tilde{T} \cap S_j} y_\gamma(\tilde{t}) > u \Bigr)\Biggr]\nonumber\\[-8pt]\\[-8pt]
&=& \frac{1}{J} \sum_{j=1}^{J} \Power_j(u),\nonumber
\end{eqnarray}
where $\Power_j(u)$ is the probability of detecting peak $j$
%
\begin{equation}
\label{eqpower-j}
\Power_j(u) = \P\Bigl\{\tilde{T} \cap S_j \ne\varnothing\mbox
{ and } \max_{t
\in\tilde{T} \cap S_j} y_\gamma(t) > u \Bigr\}.
\end{equation}
The maximum operator above indicates that if more than one significant
local maximum fall within the same peak support, only one is counted,
so power is not inflated. However, this has no effect asymptotically
because each true peak is represented by exactly one local maximum of
the smoothed observed process with probability tending to 1 (Lemma \ref
{lemmaunique-max} in Section \ref{applemmas}). The next result
indicates that both the Bonferroni and BH procedures are asymptotically
consistent. The proof is given in Section \ref{apppower}.
%
\begin{theorem}
\label{thmpower}
Let the power be defined by (\ref{eqpower}), and let $\tilde{u}_{\Bon
}$ and $\tilde{u}_{\BH}$ be the Bonferroni and BH thresholds (\ref
{eqthresh-Bon-random}) and (\ref{eqthresh-BH-random}), respectively.
Under conditions \textup{(C1)} and \textup{(C2)},
\[
\Power(\tilde{u}_{\Bon}) \to1,\qquad \Power(\tilde{u}_{\BH})
\to1.
\]
\end{theorem}

For pointwise tests, if there exists a signal anywhere, the BH
procedure is more powerful than the Bonferroni procedure
[\citet{Benjamini1995}]. This is also true in our case. Comparing
(\ref{eqthresh-Bon-fixed}) and (\ref{eqthresh-BH-fixed}), if $J \ge
1$, the threshold $u^*_{\Bon}$ is higher than the threshold
$u^*_{\BH}$, promising a larger power for the BH procedure.

\subsection{Optimal smoothing kernel}
\label{secoptimal-gamma}
The best smoothing kernel $w_\gamma(t)$ is that which maximizes the
power (\ref{eqpower}) under the true model. Because this maximization
is analytically difficult, we resort to a less formal argument here.
Lemma \ref{lemmaunique-max} in Section \ref{applemmas} shows that,
under conditions (C1) and (C2), every true peak $j$ is represented by
exactly one significant local maximum located in a small neighborhood
containing the true peak mode $\tau_j$ with probability tending to 1.
Thus the power for peak $j$ (\ref{eqpower-j}) may be approximated as
%
\begin{equation}
\label{eqapprox-power}
\Power_j(u) \approx\P\{y_\gamma(\tau_j) > u\}
= \Phi\biggl[\frac{a_j h_{j,\gamma}(\tau_j) - u}{\sigma_\gamma
}\biggr],
\end{equation}
because $y_\gamma(\tau_j) \sim N(a_j h_{j,\gamma}(\tau_j), \sigma
_\gamma^2)$. By Lemma \ref{lemmafwer-fdr-threshold} in Section \ref
{apppower}, the asymptotically equivalent thresholds (\ref
{eqthresh-Bon-fixed}) and (\ref{eqthresh-BH-fixed}) for the Bonferroni
and BH procedures satisfy $u^*_{\Bon}/a_j \to0$ and $u^*_{\BH}/a_j
\to0$ for any $j$. Thus, for large~$a_j$, the power (\ref
{eqapprox-power}) is maximized approximately by maximizing the SNR
%
\begin{equation}
\label{eqSNR}
\SNR_\gamma= \frac{a_j h_{j,\gamma}(\tau_j)}{\sigma_\gamma} =
\frac{a_j \int_{-\infty}^\infty w_\gamma(s) h_j(s) \,ds}
{\sigma\sqrt{\int_{-\infty}^{\infty} w^2_\gamma(s) \,ds}},
\end{equation}
where $\sigma$ is the standard deviation of the observed process
$y(t)$. The optimal smoothing kernel $w_\gamma(t)$ is that which is
closest to $h_j(t)$ in an $L_2$ sense. This result is similar to the
matched filter theorem for detecting a single signal peak of known
shape at a fixed time location $t$ [\citet{Pratt1991},
\citet{Simon1995}].
The result only holds approximately in our case because the peak
locations are unknown.
%
\begin{example}[(Gaussian autocorrelation model)]
\label{exGaussian-gamma-choice}
$\!\!\!$Suppose the signal peak~$j$ is a truncated Gaussian density $h_j(t) =
(1/b_j)\phi[(t-\tau_j)/b_j]\mathbf{1}[-c_j,c_j]$, $b_j,\break c_j>0$, and
let the noise be generated as in Example \ref{exGaussian-ACVF}.
Ignoring the truncation, $h_{j,\gamma}(t) = w_\gamma(t) * h_j(t)$ in
(\ref{eqSNR}) is the\vspace*{1pt} convolution of two Gaussian densities with
variances $\gamma^2$ and $b_j^2$, which is another Gaussian density
with variance $\gamma^2 + b_j^2$. Using the moments from Example \ref
{exGaussian-ACVF}, we have that
%
\begin{equation}
\label{eqSNR-Gaussian}
\SNR_\gamma= \frac{a_j h_{j,\gamma}(\tau_j)}{\sigma_\gamma} =
\frac{a_j}{\sigma\pi^{1/4}}\biggl[\frac{\gamma^2 + \nu^2}{(\gamma
^2 + b_j^2)^2}\biggr]^{1/4}.
\end{equation}
As a function of $\gamma$, the SNR is maximized at
%
\begin{equation}
\label{eqoptimal-gamma}
\argmax_\gamma\SNR_\gamma= \cases{
\sqrt{b_j^2 - 2 \nu^2}, &\quad $\nu< b_j/\sqrt{2}$, \cr
0, &\quad $\nu> b_j/\sqrt{2}$.}
\end{equation}
In particular, when $\nu=0$, we have that the optimal bandwidth for
peak $j$ is \mbox{$\gamma= b_j$}, the same as the signal bandwidth. We show
in the simulations below that the optimal $\gamma$ is indeed close to
(\ref{eqoptimal-gamma}).
\end{example}

\section{Simulation studies}
\label{secsimulations}

\subsection{Nonasymptotic performance}
\label{secsim1}
Simulations were used to evaluate the performance and limitations of
the STEM algorithm for finite range $L$ and moderate signal strength
$a$. In a segment of length $L=1000$, $J=10$ equal truncated Gaussian
peaks $a_j h_j(t) = a/b\phi[(t-\tau_j)/b] \mathbf{1}[-cb,cb]$,
$j=1,\ldots,J$, as in Example \ref{exGaussian-gamma-choice} with
$b=3$, $c=3$ and varying $a$, were placed at uniformly spaced locations
$\tau_j = (j-1/2) L/J$, $j=1,\ldots,J$, and sampled at integer values
of $t$. The noise $z(t)$ was constructed as in Example \ref
{exGaussian-ACVF} with $\sigma=1$ and varying $\nu$. Algorithm \ref
{algSTEM} was carried out using as smoothing kernel a truncated
Gaussian density $w_\gamma(t) = (1/\gamma)\phi(t/\gamma)\mathbf
{1}[-c\gamma,c\gamma]$ as in Example \ref{exGaussian-ACVF} with
$c=3$ and varying~$\gamma$. The noise parameters (\ref{eqmoments})
were estimated independently as the empirical moments of smoothed
sequences i.i.d. Gaussian noise of length 1000 and their first and
second-order differences, using the same smoothing kernel. The
Bonferroni and BH procedures were applied at level $\alpha=0.05$.


\begin{figure}

\includegraphics{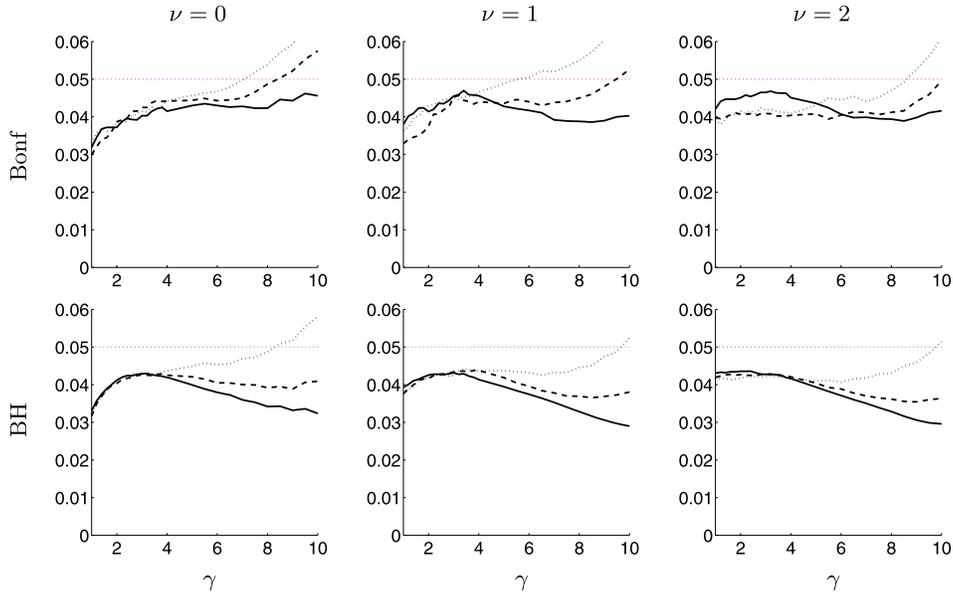}

\caption{FWER of the Bonferroni procedure (top row) and FDR of the BH
procedure (bottom row) for $a=15$ (solid), $a=12$ (dashed) and $a=9$
(dotted). Nominal error level is 0.05.}
\label{figerror}
\vspace*{-3pt}
\end{figure}

Figure \ref{figerror} shows the realized FWER and FDR levels of the
Bonferroni and BH procedures, evaluated according to (\ref{eqFWER})
and (\ref{eqFDR}) with the expectations replaced by ensemble averages
over 10,000 replications. Error rates are maintained below the nominal
level $\alpha=0.05$ for all bandwidths and large enough signal
strength $a$. The convergence is slower, however, when the bandwidth~$\gamma$
is much larger than the signal peak bandwidth $b=3$. The
increased error rates are the result of true peak maxima being shifted
from the original signal region $\mathbb{S}_1$ into the transition
region~$\mathbb{T}_{\gamma}$, where they are counted as false
positives. This phenomenon disappears with increasing signal strength~%
$a$ because the probability of obtaining any local maxima in the
transition region goes to zero asymptotically (Lemma \ref
{lemmaunique-max} in Section~\ref{applemmas}).\vadjust{\goodbreak}

\begin{figure}[b]
\vspace*{-3pt}
\includegraphics{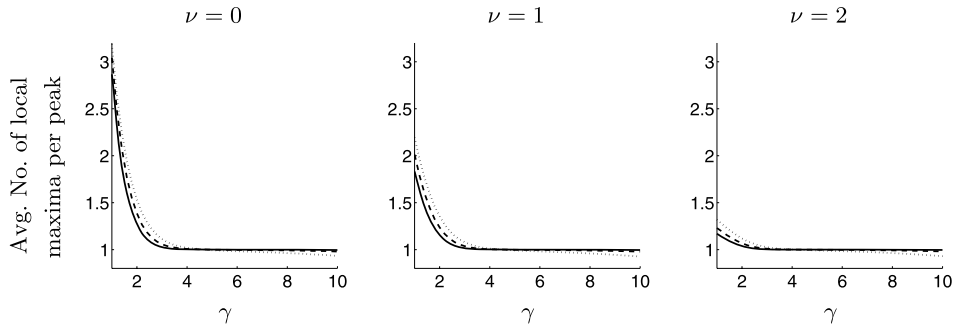}

\caption{Average number of local maxima for each true peak for $a=15$
(solid), $a=12$ (dashed) and $a=9$ (dotted).}
\label{figlocmax}
\end{figure}

As noted in Section \ref{secerrors}, each true peak may contain more
than one local maximum of the smoothed data $y_\gamma(t)$. Figure \ref
{figlocmax} shows that the expected number of local maxima per true
peak decreases with increasing bandwidth, and is essentially equal to 1
for bandwidths equal to or greater than the optimal bandwidth. It also
gets closer to 1 with increasing signal strength, consistent with the
result of Lemma \ref{lemmaunique-max}.

\begin{figure}

\includegraphics{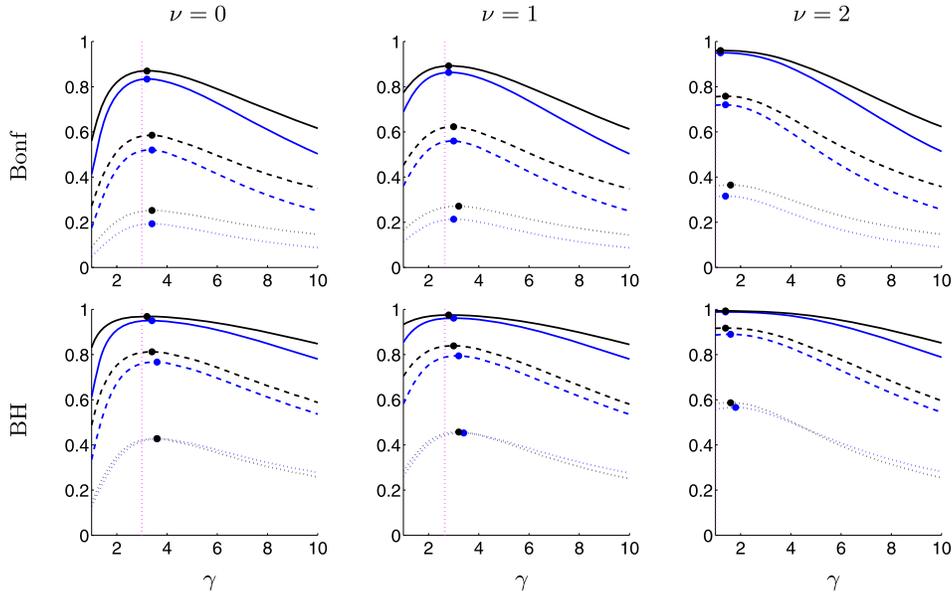}

\caption{Realized (black) and ``theoretical'' (blue)
power of the Bonferroni (top row) and BH (bottom row) procedures for
$a=15$ (solid), $a=12$ (dashed) and $a=9$ (dotted). The maxima of the
curves (solid circles) approach the asymptotic optimal bandwidth
(vertical dashed).}
\label{figpower}
\end{figure}

Figure \ref{figpower} shows the realized power of the Bonferroni and
BH procedures, evaluated according to\vadjust{\goodbreak} (\ref{eqpower}) with the
expectations replaced by ensemble averages over the same 10,000
replications. In all cases, the power increases asymptotically to~1
with the signal strength for every fixed bandwidth, and is always
larger for BH than it is for Bonferroni. The convergence is slower,
however, when the bandwidth $\gamma$ is far from the optimal value.
To understand the dependence on bandwidth, superimposed is the
theoretical approximate power (\ref{eqapprox-power}) evaluated at the
asymptotic thresholds $u^*_{\Bon}$ (\ref{eqthresh-Bon-fixed}) and
$u^*_{\BH}$ (\ref{eqthresh-BH-fixed}) and plugging in the SNR (\ref
{eqSNR-Gaussian}). The ``theoretical'' power curves largely capture the
shape of the realized ones, but are lower because the asymptotic
thresholds are more conservative. The curve shape is mostly determined
by the SNR (\ref{eqSNR-Gaussian}) as a function of $\gamma$. The
bandwidth $\gamma$ producing the largest power is always larger than
the theoretical optimal bandwidth (\ref{eqoptimal-gamma}), but it
approaches it from the right as $a$ increases.

%
\begin{figure}[b]

\includegraphics{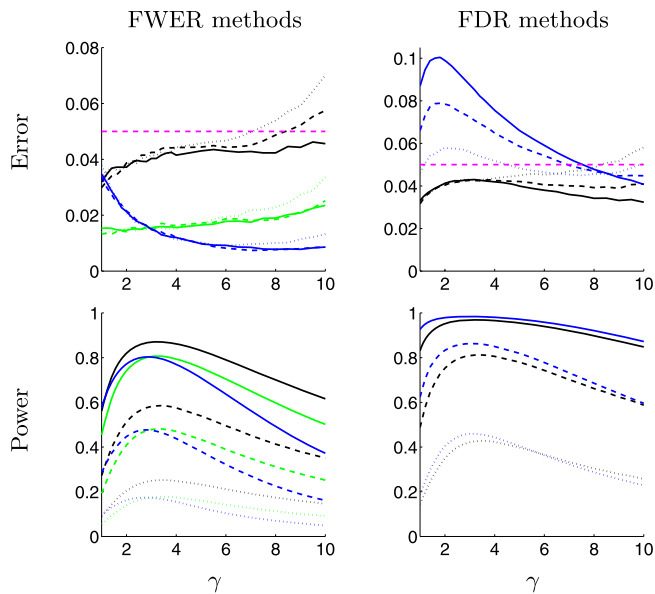}

\caption{Left panels: FWER and power of three FWER methods: STEM with
Bonferroni (black), Bonferroni on all $L$ locations (blue) and Supremum
(green). Right panels: FDR and power of three FDR methods: STEM with BH
(black), BH on all $L$ locations (blue). Results in all panels are for
$a=15$ (solid), $a=12$ (dashed) and $a=9$ (dotted). Nominal error level
is 0.05.}
\label{figcompare-methods}
\end{figure}
%

%

\subsection{Unequal peaks}
By assumption (Section \ref{secmodel}), the signal peaks need not be
equal. As in Figure \ref{figsimulexample}, $J=5$ unequal peaks
(Epanechnikov, triangular and truncated Gaussian, Laplace and Cauchy,
with average half-support 24) were corrupted with white standard normal
noise. Algorithm \ref{algSTEM} was applied using a quartic smoothing
kernel $w_\gamma(t) = 15/(16\gamma)[1 - (t/\gamma)^2]^2 \mathbf
{1}[-\gamma,\gamma]$ with varying $\gamma$, the noise parameters
estimated independently as in Section \ref{secsim1}. For this
configuration and 10,000 repetitions, the error was controlled below
the nominal level 0.05 for values of $\gamma$ up to 40, obtaining a
maximum power of 0.81 and 0.88 for the Bonferroni and BH procedures at
$\gamma=18$. The maximizing bandwidth represents the average best
match between the quartic smoothing kernel and the peaks present in the data.

\subsection{Overlapping peaks}
The theory of Section \ref{sectheory} assumed that the signal peaks
had nonoverlapping supports. Simulations similar to those of Section~%
\ref{secsim1} with $J=10$ partially overlapping peaks showed that the
error rates were below the nominal level regardless of the amount of
overlap between peaks. The detection power, however, deceptively
increased with increasing overlap. This is because definition (\ref
{eqpower}) counts two overlapping peaks as detected even if only one
significant local maximum is found in the overlapping region between
them, as it belongs to both. Definition (\ref{eqpower}) does not
measure the ability to distinguish between overlapping peaks.

\subsection{Comparison with pointwise testing}
\label{seccompare-simul} To see the benefits of testing local maxima,
Figure \ref{figcompare-methods} compares the performance of the STEM
algorithm (with Bonferroni and BH corrections) to three other methods
that test at every single location. Simulated data sets as in Section
\ref{secsim1} with $b=3$ and $\nu= 0$ were smoothed with varying
$\gamma$. For the pointwise Bonferroni and BH methods, $p$-values for
testing $H_0\dvtx \mu(t) = 0$ at each $t=1,\ldots,L=1000$ were computed as
$p(t) = 1-\Phi[y_\gamma(t)/\sigma_\gamma]$ and then corrected using
Bonferroni and BH, respectively. The method ``Supremum'' was adapted
from \citet{Worsley1996b} as follows. The probability that the
supremum of any differentiable random process $f(t)$ in the interval
$[0,T]$ exceeds $u$ is bounded by [\citet{Adler2007}]
%
\begin{equation}
\label{eqsup}
\P\Bigl(\sup_{t\in[0,T]}f(t) \ge u\Bigr) \le\P[f(0) \ge
u] + \E[N_u],
\end{equation}
where $N_u$ is the number of up-crossings by $f(t)$ of the level $u$ in
$[0, T]$. For the stationary Gaussian process $z_\gamma(t)$,
application of the Kac--Rice formula [\citet{Cramer1967}, page 194]
gives that $\E[N_u] = L (\sqrt{\lambda_{2, \gamma}}/\allowbreak\sigma_\gamma)
\phi(u/\sigma_\gamma)$. The significance threshold is found as the
largest $u$ such that
%
\begin{equation}
\label{eqsup-approx}
\P\Bigl(\sup_{t\in[-L/2,L/2]} z_\gamma(t) \ge u\Bigr) \le1
-\Phi\biggl(\frac{u}{\sigma_\gamma}\biggr) +
L \frac{\sqrt{\lambda_{2,\gamma}}}{\sigma_\gamma} \phi
\biggl(\frac{u}{\sigma_\gamma}\biggr) \le\alpha.
\end{equation}

Figure \ref{figcompare-methods} indicates that the pointwise
Bonferroni correction is too conservative. The Supremum method, despite
accounting explicitly for the noise autocorrelation, performs only
slightly better than pointwise Bonferroni, and not as well as
Bonferroni performed on local maxima. The pointwise BH correction is
designed to control FDR at the level of individual locations, and thus
produces too many false positives when the FDR is measured in terms of
detected peaks using (\ref{eqFDR}). Further simulations with $\nu=1$
and $\nu=2$ yielded similar results (not shown).

\subsection{Automatic bandwidth selection}
\begin{figure}
\begin{tabular}{@{}c@{\quad}c@{}}

\includegraphics{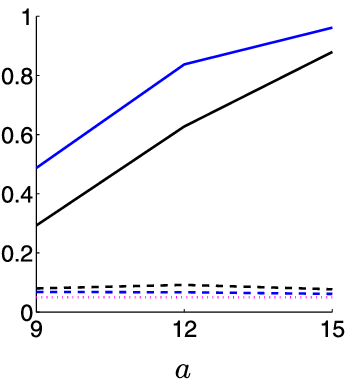}
 & \includegraphics{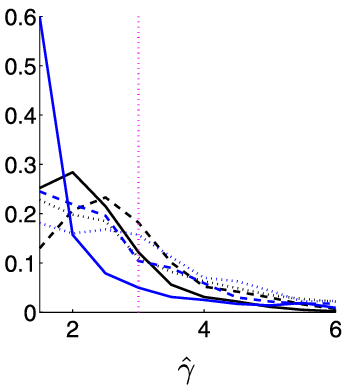}\\
(a) & (b)
\end{tabular}
\caption{\textup{(a)} Power (solid) and realized error rate (dashed)
for Bonferroni (black) and BH (blue) with automatic bandwidth selection
as a function of signal strength $a$. Nominal error level is 0.05.
\textup{(b)} Proportion of automatically chosen smoothing bandwidth
$\hat{\gamma}$ over 1000 simulations for Bonferroni (black) and BH
(blue); results are for $a=15$ (solid), $a=12$ (dashed) and $a=9$
(dotted). Nominal optimal bandwidth is $\gamma=3$.} \label{figgammaest}
\end{figure}

Rather than using a fixed smoothing bandwidth $\gamma$, the bandwidth
may be chosen automatically from the data as the one that yields the
largest number of discoveries for a fixed error level. For simulated
data sets as in Section \ref{secsim1} with $b=3$ and $\nu= 0$, the
STEM algorithm was applied with $\gamma$ ranging from $\gamma=b/2=1.5$ to
$\gamma=2b=6$, and results were retained for the bandwidth $\hat{\gamma}$
that yielded the largest number of discoveries in each run.
Figure \ref{figgammaest}(a) shows that this automatic criterion biases
the results toward more detected peaks and therefore results in higher
%
error rates (and power) than those obtained when $\gamma$ is fixed
(Figure \ref{figpower}). It also tends to select bandwidths that are
smaller than the nominal optimal value $\gamma=b$ [Figure \ref
{figgammaest}(b)], with averages ranging between about 2.1 and 2.9.

\section{Data example}
\label{secdata} 
\begin{figure}

\includegraphics{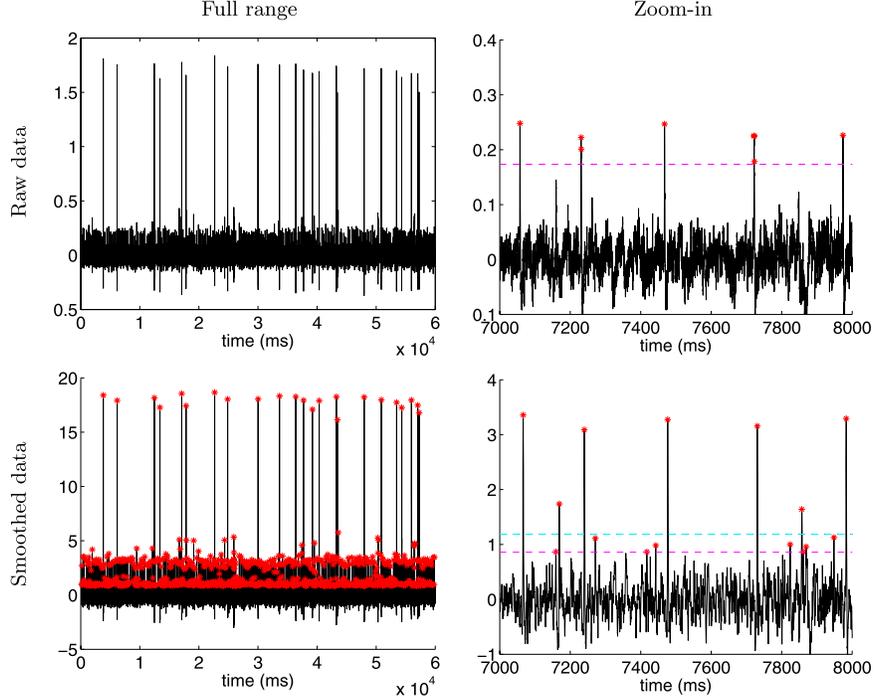}
 \caption{Top row: the neural spike data (test set);
the stars in the right panel indicate peaks that are higher than 4
standard deviations of the raw data (dashed line), as suggested by
Segev et al. (\protect\citeyear{Segev2004}). Bottom
row: the data smoothed using the estimated peak shape as kernel; the
stars indicate significant local maxima higher than the BH threshold
(magenta dashed line) at level 0.01. The Bonferroni threshold is
indicated by the cyan dashed line.} \label{figspike-data}
\end{figure}

The data consists of recordings from a single electrode
inserted in a salamander's retina, digitized at a sampling frequency of
10 kHz. Data of these kind are routinely collected in large amounts in
neuroscience experiments [\citet{Baccus2002},
Segev et al. (\citeyear{Segev2004})]. For the purposes of this paper, three data
sets were used:
\begin{longlist}[(3)]
\item[(1)] Test set: 60 seconds of recordings of live cells in the
dark.
\item[(2)] Training set 1: 60 seconds of recordings of live
cells in the dark.
\item[(3)] Training set 2: 60 seconds of recordings
after the retina was allowed to die.
\end{longlist}
Each period of 60 seconds corresponds to $L = 6\times10^5$ samples.
The goal of the analysis was to detect neuronal spikes in the test set
(Figure \ref{figspike-data}, top left).

\begin{figure}
\begin{tabular}{@{}c@{\quad}c@{}}

\includegraphics{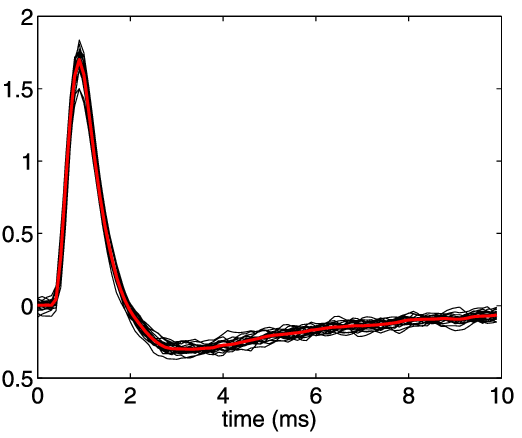}
 & \includegraphics{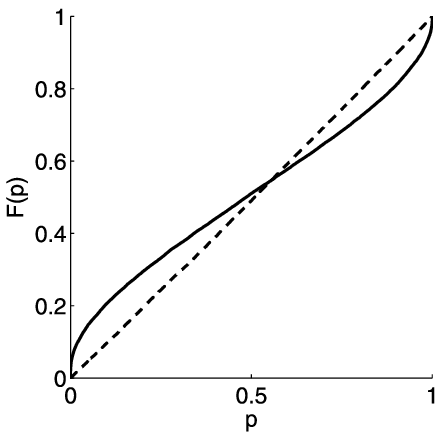}\\
(a) & (b)
\end{tabular}
\caption{\textup{(a)} 23 strongest spikes aligned by their maximum (black);
their average (red) is the estimated template. \textup{(b)} Empirical cdf of
$p$-values for the test set (solid) and training set 2 (dashed).}
\label{figtemplate-Fp}
\end{figure}

Assuming that neuronal action potentials have similar shapes, to
maximize the SNR (\ref{eqSNR}), the smoothing filter should be close
in shape and bandwidth to that of the peaks to be detected. Training
set 1 was used to estimate the peak shape. In training set 1, spikes
with raw maximum exceeding 1 were selected and aligned by their maxima
[Figure \ref{figtemplate-Fp}(a)]. The peak shape template was obtained
as the average of the 23 selected major spikes and truncated to a
length of 100 samples.

Training set 2, recorded under pure noise conditions, was used to
estimate the noise parameters. The noise in training set 2 can be well
modeled by an AR(3) process with autoregressive coefficients $-$1.13,
0.42 and $-$0.13, estimated by the Yule--Walker algorithm, so that
whitening with these coefficients produces a process whose
autocovariance function cannot be distinguished from that of white
noise using a Bartlett's test. A similar analysis in segments of length
$L/10$ showed that the estimated AR coefficients have a coefficient of
variation of no more than 1\% over the 10 segments, supporting the
stationarity assumption. A Jarque--Bera test of normality for the
entire sequence returned a $p$-value of 0.224, supporting the Gaussianity
assumption.

Convolving\vspace*{1pt} training set 2 with the template of Figure \ref
{figtemplate-Fp}(a) produced smoothed noise with spectral moments $\hat
{\sigma}_\gamma^2 = 4.22\times10^{-4}$, $\hat{\lambda}_{2,\gamma }^2 =
1.20\times10^{-4}$ and $\hat{\lambda}_{4,\gamma}^2 =
1.96\times10^{-4}$, estimated respectively by the empirical variances
of the observed process, its first-order difference and its
second-order difference. Given the length of the process, the standard
error
of these estimates is negligible.

Algorithm \ref{algSTEM} was applied to the test set (Figure \ref
{figspike-data}, top left) by convolving it with the template of Figure
\ref{figtemplate-Fp}(a), producing the smoothed process in Figure~\ref
{figspike-data} (bottom left). In $L = 6\times10^5$ samples, $\tilde
{m} = 30\mbox{,}426$ local maxima were found and their $p$-values were computed
according to (\ref{eqp-value}) and (\ref{eqdistr}), plugging in the
estimates $\hat{\sigma}_\gamma^2$, $\hat{\lambda}_{2,\gamma}^2$
and $\hat{\lambda}_{4,\gamma}$ found above. The empirical\vspace*{1pt} cdf of the
$p$-values [Figure \ref{figtemplate-Fp}(b)] shows a large fraction of
nonnull $p$-values near 0. For comparison, the same procedure of
smoothing, finding local maxima and computing their $p$-values was
applied to training set 2. The empirical cdf of those $p$-values is
virtually uniform, emphasizing that formula (\ref{eqdistr}) for
Gaussian noise is appropriate. Also in Figure \ref{figtemplate-Fp}(b),
the excess of large $p$-values from the test set is due to the negative
portions of the smoothing function [Figure \ref{figtemplate-Fp}(a)].
These produce small negative anti-spikes whose $p$-values are large when
tested for positiveness.

Applying the BH procedure to the $\tilde{m} = 30\mbox{,}426$ $p$-values obtained
from the test set at FDR level 0.01 resulted in a $p$-value threshold of
$2.76\times10^{-4}$ and $R = 843$ significant local maxima. These are
indicated in Figure \ref{figspike-data} (bottom left), showing three
levels of spike strengths. Figure \ref{figspike-data} (bottom right)
zooms in to show a few of the weaker spikes. Applying the Bonferroni
procedure instead in Algorithm \ref{algSTEM} resulted in a $p$-value threshold of
$3.29\times10^{-7}$ and only 411 detected spikes.

For comparison, Figure \ref{figspike-data} (top right) shows the same
segment of the raw data and the spikes selected using one of the
recommended methods in the neuroscience literature, which is to
threshold at 4 standard deviations of the raw data [Segev et al.
(\citeyear{Segev2004})]. Our method is able to identify more
spikes at a low FDR level of 0.01, but more importantly, it attaches to
the findings a~significance level, expecting about 1\% of the detected
spikes to be false. The conventional method does not offer this useful
statistical interpretation.

As in Section \ref{seccompare-simul}, computing $p$-values at each
location as $p(t) = 1 - \Phi(y_\gamma(t)/\hat{\sigma}_\gamma)$,
$t=1,\ldots,L$, and applying a global Bonferroni at level 0.01 was
more conservative, resulting in a height threshold of 1.235 (comparable
to Figure~\ref{figspike-data} bottom right) and detecting only 393
spikes. Similarly, the ``Supremum'' method, applied by replacing $\hat
{\sigma}_\gamma$ and $\hat{\lambda}_\gamma$ in (\ref
{eqsup-approx}) at level 0.01, yielded a height threshold 1.229 and 394
detected spikes. Finally, applying the global BH procedure at level
0.01 with $L$ $p$-values gave a height threshold of 0.780 detecting 1149
spikes, but as shown in Section \ref{seccompare-simul}, this result is
too optimistic because the actual error rate for peaks is higher than 0.01.

\section{Discussion}
\label{secdiscussion}


For the theoretical results, the most critical assumptions were that
the noise process is stationary ergodic Gaussian and that the signal
peaks are unimodal with compact support. The Gaussianity assumption was
chosen because it enabled a closed formula for computing the $p$-values
associated with the heights of local maxima. For non-Gaussian noise,
$p$-values could be computed via Monte Carlo simulation.

The assumption of compact support for the signal peaks was necessary
for true and false positives to be well defined. \citet{Chumbley2010}
argued for testing local maxima when the signal spreads over the entire
domain, but in that case every positive is a true positive, making the
inference unclear. On the other hand, agreeing with \citet
{Chumbley2009}, applying BH globally resulted in inflated error rates
for peaks, while applying Bonferroni or the Supremum method globally
was too conservative. The unimodality assumption made local maxima good
representatives of true peaks, being unique for medium to large
bandwidths and asymptotically for increasing signal strength.

The strong signal assumption in condition (C2) was introduced to remove
the excess error produced by the smoothed signal spreading into the
neighboring null regions, thereby enabling asymptotic error control.
The assumption is not restrictive in the sense that the search space
may grow exponentially faster. Similar conditions are common for
high-dimensional data. If the data are pointwise test statistics based
on a sample of size $n$, with SNR increasing as $a = \sqrt{n} \to
\infty$, then the condition $(\log L)/a^2 \to0$ becomes $(\log L)/n
\to0$. This is similar to the condition $(\log p)/n \to0$ required for
consistent model selection in high-dimensional regression under
sparsity where $p$ is the number of features [\citet{Candes2007},
\citet{Zhang2010}]. Our results, however, do not require sparsity.
Condition (C2) is easy to state but stronger than needed; upon close
inspection of the proof of Lemma \ref{lemmaunique-max} in Section
\ref{applemmas}, the limit of $(\log L)/a^2$ need not be zero but need
only be bounded by a constant that depends on the signal and noise
first and second derivatives.

While the theory was developed for continuous processes, in practice
the observations are given in a discrete grid. In our simulations we
found that the results were not reliable when the smoothing bandwidth
was smaller than the grid spacing, as the theory for continuous random
processes is no longer a good approximation in that case.

The asymptotic error control and power consistency did not require the
peaks to have the same shape or width. The asymptotic results were
found to hold in practice for a wide range of bandwithds and strong
enough signal. However, the convergence rate was slower for bandwidths
less than half or more than double the optimal value. The matched
filter principle suggests that the smoothing kernel should be chosen to
be as close as possible in an $L^2$ sense to the peaks to be detected.
In the neuronal data analyzed, the peak shape and width were estimated
from the data, dictating the best smoothing kernel. If the peaks to be
detected have different widths, then the bandwidth may be adapted to
the width of each peak. We leave this possibility for future work, as
well as the obliged extension of the proposed methods to two- and
three-dimensional domains.

\section{Technical details}
\label{secproofs}


\subsection{Supporting results}
\label{applemmas}

%
\begin{lemma}
\label{lemmaPalm-S0}
Let $\tilde{m}_{0,\gamma} = \#\{ t\in\tilde{T}\cap\mathbb
{S}_{0,\gamma}\}$ be be the number of local maxima of $y_\gamma(t)$
[or $z_\gamma(t)$] in $\mathbb{S}_{0,\gamma}$. Let $V_\gamma(u) = \#
\{t\in\tilde{T}\cap\mathbb{S}_{0, \gamma} \dvtx y_\gamma(t) > u\}$ be
the number of local maxima of $y_\gamma(t)$ [or $z_\gamma(t)$] in
$\mathbb{S}_{0,\gamma}$ whose heights are above the level $u$. Then
\[
\frac{V_\gamma(u)}{\tilde{m}_{0,\gamma}} \to\frac{\E[V_\gamma
(u)]}{\E[\tilde{m}_{0,\gamma}]} = F_\gamma(u)
\]
in probability as $L\to\infty$, where $F_\gamma(u)$ is the Palm
distribution (\ref{eqpalm}).
\end{lemma}
\begin{pf}
Notice that $y_\gamma(t) = z_\gamma(t)$ for all $t \in\mathbb
{S}_{0, \gamma}$, so the process $y_\gamma(t)$ has the same
properties as the stationary process $z_\gamma(t)$ on the set $\mathbb
{S}_{0, \gamma}$. By ergodicity, the weak law of large numbers applied
to the numerator and denominator gives that
%
\begin{equation}
\label{eqV-over-m}
\frac{V_\gamma(u)}{\tilde{m}_{0,\gamma}} =
\frac{\#\{t\in\tilde{T}\cap\mathbb{S}_{0, \gamma} \dvtx z_\gamma(t)
> u\}/L}{\#\{ t\in\tilde{T}\cap\mathbb{S}_{0,\gamma}\}/L}
\end{equation}
converges to [\citet{Cramer1967}]
\[
\frac{\E[\#\{t\in\tilde{T}\cap\mathbb{S}_{0, \gamma} \dvtx z_\gamma
(t) > u\}]}{\E[\#\{ t\in\tilde{T}\cap\mathbb{S}_{0,\gamma}\}]} =
\frac{\E[V_\gamma(u)]}{\E[\tilde{m}_{0,\gamma}]}.
\]
But also by ergodicity, ratio (\ref{eqV-over-m}) converges to the
conditional probability $\P[z_\gamma(t) > u | t\in\tilde{T}\cap
\mathbb{S}_{0,\gamma}] = F_\gamma(u)$ by Definition (\ref{eqpalm}).
The two limits must be equal.
\end{pf}
%
%
\begin{lemma}
\label{lemmabounds}
Assume the model of Section \ref{secmodel}. Let $S_{j,\gamma} =
I_j^{\mathrm{left}} \cup I_j^{\mathrm{mode}} \cup I_j^{\mathrm{right}}$ be a
partition, where $I_j^{\mathrm{mode}} = [c_j,d_j] \subset S_j$ is a fixed
interval containing the mode of $\mu_\gamma(t) = a_j h_{j,\gamma
}(t)$ in $S_j$ as an interior point, such that $\ddot{h}_{j,\gamma
}(t) < 0$ for $t \in I_j^{\mathrm{mode}}$, $\dot{h}_{j,\gamma}(t) > 0$ for
$t \in I_j^{\mathrm{left}}$ and $\dot{h}_{j,\gamma}(t) < 0$ for $t \in
I_j^{\mathrm{right}}$. Let:
\begin{itemize}
\item$M_j$ be the largest value of $|h_{j,\gamma}(t)|$ in
$S_{j,\gamma}$;
\item$C_j$ be the smallest value of $|\dot{h}_{j,\gamma}(t)|$ in
$I_j^{\mathrm{side}} = I_j^{\mathrm{left}} \cup I_j^{\mathrm{right}}$;
\item$D_j$ be the smallest value of $|\ddot{h}_{j,\gamma}(t)|$
in\vadjust{\goodbreak}
$I_j^{\mathrm{mode}}$.
\end{itemize}
For $\tilde{T}$ given by (\ref{eqT}) and any threshold $u$,
%
\begin{eqnarray}
\label{eqI-bound}
&&
\P(\# \{t \in\tilde{T}\cap I_j^{\mathrm{side}}\} = 0) \nonumber\\
&&\qquad\ge
2\Phi\biggl(\frac{a_j C_j}{\sqrt{\lambda_{2,\gamma}}}\biggr) -1 -
|I_j^{\mathrm{side}}| \sqrt{\frac{\lambda_{4, \gamma} }{
\lambda_{2, \gamma}}} \phi\biggl(\frac{a_j C_j}{\sqrt{\lambda_{2,
\gamma}}}\biggr), \nonumber\\
&&\P(\# \{t \in\tilde{T}\cap I_j^{\mathrm{mode}}\} =
1)\nonumber\\[-8pt]\\[-8pt]
&&\qquad\ge
\Phi\biggl(\frac{a_j D_j}{\sqrt{\lambda_{4,\gamma}}}\biggr) -
|I_j^{\mathrm{mode}}| \sqrt{\frac{\lambda_{6, \gamma} }{
\lambda_{4, \gamma}}} \phi\biggl(\frac{a_j D_j}{\sqrt{\lambda_{4,
\gamma}}}\biggr) -
2\Phi\biggl(\frac{-a_j C_j}{\sqrt{\lambda_{2,\gamma}}}\biggr), \nonumber\\
&&\P\bigl(\# \{t \in\tilde{T}\cap I_j^{\mathrm{mode}}\dvtx y_\gamma(t) > u \}
= 1\bigr) \nonumber\\
&&\qquad\ge
\Phi\biggl(\frac{a_j D_j}{\sqrt{\lambda_{4,\gamma}}}\biggr) -
|I_j^{\mathrm{mode}}| \sqrt{\frac{\lambda_{6, \gamma} }{
\lambda_{4, \gamma}}} \phi\biggl(\frac{a_j D_j}{\sqrt{\lambda_{4,
\gamma}}}\biggr) -
\Phi\biggl(\frac{u-a_j M_j}{\sigma_\gamma}\biggr),
\nonumber
\end{eqnarray}
where $\sigma_\gamma$, $\lambda_{2, \gamma}$ and $\lambda_{4,
\gamma}$
are given by (\ref{eqmoments}) and $\lambda_{6, \gamma} =
\E[\dddot_\gamma(t)]$.
\end{lemma}
\begin{pf}
(1) Consider first the compact interval $I_j^{\mathrm{left}}$. The
probability that there are no local maxima of $y_\gamma(t)$ in
$I_j^{\mathrm{left}}$ is greater\vspace*{1pt} than the probability that $\dot{y}_\gamma
(t) > 0$ for all $t$ in the interval. This probability is equal to
%
\begin{eqnarray}
\label{eqI-left}
\P\Bigl(\inf_{I_j^{\mathrm{left}}} \dot{y}_\gamma(t) > 0\Bigr)
&\ge&
\P\Bigl(\inf_{I_j^{\mathrm{left}}} \dot{z}_\gamma(t) > -\inf_{I_j^{\mathrm{left}}} \dot{\mu}_\gamma(t) \Bigr)
\nonumber\\[-8pt]\\[-8pt]
&=&
1 - \P\Bigl(\sup_{I_j^{\mathrm{left}}} [-\dot{z}_\gamma(t)] > a_j
C_j^{\mathrm{left}} \Bigr),\nonumber
\end{eqnarray}
where $C_j^{\mathrm{left}} > 0$ is the smallest value of
$\dot{h}_{j,\gamma}(t)$ in $I_j^{\mathrm{left}}$. Inequality (\ref{eqsup})
applies above to the stationary Gaussian process $-\dot{z}_\gamma(t)$.
The Kac--Rice formula [\citet{Cramer1967}, page 194] gives in this case
that $\E[N_u] = |I_j^{\mathrm{left}}| \sqrt{\lambda_{4,
\gamma}}/\sqrt{\lambda_{2, \gamma}} \phi(u/\sqrt{\lambda_{2,
\gamma}})$. Thus (\ref{eqI-left}) has the lower bound
\[
\P(\# \{t \in\tilde{T}\cap I_j^{\mathrm{left}}\} = 0) \ge
\Phi\biggl(\frac{a_j C_j^{\mathrm{left}}}{\sqrt{\lambda_{2,\gamma
}}}\biggr) -
|I_j^{\mathrm{left}}| \sqrt{\frac{\lambda_{4, \gamma} }{
\lambda_{2, \gamma}} } \phi\biggl(\frac{a_j C_j^{\mathrm{left}}}{\sqrt
{\lambda_{2, \gamma}}}\biggr).
\]
A similar\vspace*{1pt} calculation for $I_j^{\mathrm{right}}$ gives a similar bound with
the superscript ``left'' replaced by ``right'' and $C_j^{\mathrm{right}} >
0$ being the smallest value of $|\dot{h}_{j,\gamma}(t)|$ in $I_j^{\mathrm{right}}$. Putting the two together, the required probability $\P(\# \{t
\in\tilde{T}\cap I_j^{\mathrm{side}}\})$ that there are no local maxima in
$I_j^{\mathrm{left}}$ nor $I_j^{\mathrm{right}}$ is bounded as in the first row
of (\ref{eqI-bound}).

(2) The probability\vspace*{1pt} that $y_\gamma(t)$ has no local maxima in
$I_j^{\mathrm{mode}}$ is less than the probability that $\dot{y}_\gamma
(c_j) \le0$ or $\dot{y}_\gamma(d_j) \ge0$, for a positive
derivative at $c_j$ and a~negative one at $d_j$ would imply the
existence of at least one local maximum in $I_j$. Thus, the probability
of no local maxima in $I_j^{\mathrm{mode}}$ is bounded above~by
%
\begin{eqnarray}
\label{eqI-mode-bound-0}\quad
\P(\#\{t \in\tilde{T}\cap I_j^{\mathrm{mode}}\} = 0)
&\le&\P[\dot{y}_\gamma(c_{j}) \le0] +
\P[\dot{y}_\gamma(d_{j}) \ge0] \nonumber\\
&=& \Phi\biggl(\frac{-a_j \dot{h}_{j,\gamma}(c_j)}{\sqrt{\lambda
_{2,\gamma}}}\biggr) +
1 - \Phi\biggl(\frac{-a_j \dot{h}_{j,\gamma}(d_j)}{\sqrt{\lambda
_{2,\gamma}}}\biggr)\\
&\le&2 - 2\Phi\biggl(\frac{a_j C_j}{\sqrt{\lambda_{2,\gamma
}}}\biggr),
\nonumber
\end{eqnarray}
because $\dot{y}_\gamma(t) \sim N(\dot{\mu}_\gamma(t), \lambda
_{2,\gamma})$ and $\dot{h}_\gamma(c_j) > C_j > 0$ and $-\dot
{h}_\gamma(d_j) > C_j > 0$.

On the other hand, the probability that $y_\gamma(t)$ has more than one
local maxima in $I_j^{\mathrm{mode}}$ is less\vspace*{2pt} than the
probability that $\ddot{y}_\gamma(t)> 0$ for some $t$ in~$I_j^{\mathrm{mode}}$. This probability is
\[
\P\Bigl(\sup_{I_j^{\mathrm{mode}}} \ddot{y}_\gamma(t) > 0\Bigr) \le
\P\Bigl(\sup_{I_j^{\mathrm{mode}}} \ddot{z}_\gamma(t) > a_j D_j \Bigr),
\]
where $D_j < 0$ is the largest value of $\ddot{\mu}_\gamma(t) < 0$
in $I_j^{\mathrm{mode}}$. Applying (\ref{eqsup}) to the process $\ddot
{z}_\gamma(t)$ gives the further upper bound
%
\begin{equation}
\label{eqI-mode-bound-1}\qquad
\P(\# \{t \in\tilde{T}\cap I_j^{\mathrm{mode}}\} \ge1)
\le
1 - \Phi\biggl(\frac{a_j D_j}{\sqrt{\lambda_{4,\gamma}}}\biggr) +
|I_j^{\mathrm{mode}}| \sqrt{\frac{\lambda_{6, \gamma} }{
\lambda_{4, \gamma}}} \phi\biggl(\frac{a_j D_j}{\sqrt{\lambda_{4,
\gamma}}}\biggr).
\end{equation}
Putting (\ref{eqI-mode-bound-0}) and (\ref{eqI-mode-bound-1})
together gives the bound in the second row of (\ref{eqI-bound}).

(3) The probability that no local maxima of $y_\gamma(t)$ in $I_j^{\mathrm{mode}}$ exceed the threshold $u$ is less than the probability that
$y_\gamma(t)$ is below $u$ anywhere in~$I_j^{\mathrm{mode}}$, so it is
bounded above by $\Phi[(u-a_j M_j)/\sigma_\gamma]$. On the other
hand, the probability that more than one local maxima of $y_\gamma(t)$
in $I_j^{\mathrm{mode}}$ exceed\vspace*{1pt} $u$ is less than the probability that there
exist more than one local maximum, which is bounded above by (\ref
{eqI-mode-bound-1}). Putting the two together gives the bound in the
third row of (\ref{eqI-bound}).
\end{pf}
%
%
\begin{lemma}
\label{lemmaunique-max} Assume the model of Section \ref{secmodel}. For
$\tilde{T}$ given by (\ref{eqT}), let $\tilde{m}_{1, \gamma} = \#\{
\tilde{T}\cap \mathbb{S}_{1, \gamma} \}$ be the number of local maxima
in the set $\mathbb{S}_{1, \gamma}$, and recall that $W_\gamma(u) =
\#\{t\in \tilde{T}\cap\mathbb{S}_{1, \gamma}\dvtx y_\gamma( t)> u \}$
is the number of local maxima in~$\mathbb{S}_{1, \gamma}$ above
threshold~$u$. Under conditions \textup{(C1)} and \textup{(C2)}:
\begin{longlist}[(5)]
\item[(1)]
The probability that $y_\gamma(t)$ has any local maxima in the
transition region $\mathbb{T}_\gamma$ tends to 0.
\[
\P(\# \{t\in\tilde{T}\cap\mathbb{T}_\gamma\} \ge1 )
\to0.
\]
\item[(2)]
The probability to get exactly $J$ local maxima in the
set $\mathbb{S}_{1, \gamma}$,
\[
\P( \tilde{m}_{1, \gamma} = J) = \P(\#\{t \in
\tilde{T}\cap\mathbb{S}_{1, \gamma} \} = J) \to1.\vadjust{\goodbreak}
\]
\item[(3)] The probability to get exactly $J$ local maxima in the
set $\mathbb{S}_{1, \gamma}$ that exceed any fixed threshold $u$,
\[
\P[ W_\gamma(u) = J] = \P[\# \{t\in\tilde
{T}\cap\mathbb{S}_{1, \gamma}\dvtx y_\gamma(t)> u \}=J ] \to1.
\]
\item[(4)] $\tilde{m}_{1, \gamma}/L \to A_1$ in probability.
\item[(5)] $W_\gamma(u)/\tilde{m}_{1, \gamma} \to1$ in probability.
\end{longlist}
\end{lemma}
\begin{pf}
(1)
Write $\mathbb{T}_\gamma= \bigcup_{j=1}^{J} T_{j,\gamma}$, where
$T_{j,\gamma} = S_{j,\gamma} \setminus S_j$ is the transition region
for peak $j$ (Figure \ref{figs0}). Under the assumptions of Lemma \ref
{lemmabounds}, $T_{j,\gamma}$ is a~subset of $I_j^{\mathrm{side}}$ because
$I_j^{\mathrm{left}}$ or $I_j^{\mathrm{right}}$ may include points inside $S_j$.
Using~(\ref{eqI-bound}), the required probability $\P(\# \{t\in
\tilde{T}\cap\mathbb{T}_\gamma\} \ge1 )$ that $y_\gamma
(t)$ has any local maxima in the transition region $\mathbb{T}_\gamma
$ is bounded above by
\begin{eqnarray*}
&&\sum_{j=1}^{J} \biggl[
2 - 2\Phi\biggl(\frac{a_j C_j}{\sqrt{\lambda_{2,\gamma}}}\biggr) +
|I_j^{\mathrm{side}}| \sqrt{\frac{\lambda_{4, \gamma} }{
\lambda_{2, \gamma}}} \phi\biggl(\frac{a_j C_j}{\sqrt{\lambda_{2,
\gamma}}}\biggr) \biggr] \\
&&\qquad\le2 \frac{J}{L} L \biggl[1 - \Phi\biggl(\frac{a C}{\sqrt{\lambda
_{2,\gamma}}}\biggr) \biggr] +
L \sqrt{\frac{\lambda_{4, \gamma} }{ \lambda_{2, \gamma}}} \phi
\biggl(\frac{a C}{\sqrt{\lambda_{2, \gamma}}}\biggr),
\end{eqnarray*}
where $a > 0$ is the infimum of the $a_j$'s and $C > 0$ is the infimum
of the~$C_j$'s, that is, the infimum of $|\dot{h}_{j,\gamma}(t)|$ for
$t \in\bigcup_{j=1}^{J} I_j^{\mathrm{side}}$ [recall that every peak~%
$h_{j,\gamma}(t)$ has no critical points in the transition region for
any $j$]. But the expression above goes to zero under condition (C2)
because, for any $K > 0$,
\[
L\phi(Ka) = \frac{1}{\sqrt{2\pi}} \exp\biggl[ a^2\biggl( \frac
{\log L}{a^2} - \frac{K^2}{2}\biggr)\biggr] \to0
\]
and $L[1-\Phi(Ka)] \le L\phi(Ka)/(Ka) \to0$.

(2) The required probability to obtain exactly $J$ local maxima in the
set $\mathbb{S}_{1, \gamma} = \bigcup_{j=1}^{J} S_{j, \gamma}$ is
greater than the probability of obtaining exactly one local maximum in
each interval $I_j^{\mathrm{mode}} \subset S_j$ and
none in $I_j^{\mathrm{side}}$ for any $j$. Thus,
using~(\ref{eqI-bound}), the required probability is bounded below by
\begin{eqnarray*}
&& \P\Biggl[\bigcap_{j=1}^{J} (\#\{t \in\tilde{T}\cap I_j^{\mathrm{mode}} \} = 1
\cap\#\{t \in\tilde{T}\cap I_j^{\mathrm{side}} \} = 0 )
\Biggr] \\
&&\qquad\ge1- \sum_{j=1}^{J} [1 - \P(\#\{t \in\tilde{T}\cap
I_j^{\mathrm{mode}} \} = 1 \cap\#\{t \in\tilde{T}\cap I_j^{\mathrm{side}}
\} = 0 ) ] \\
&&\qquad\ge1- \sum_{j=1}^{J} \Biggl[
5 - 4\Phi\biggl(\frac{a_j C_j}{\sqrt{\lambda_{2,\gamma}}}\biggr)
- \Phi\biggl(\frac{a_j D_j}{\sqrt{\lambda_{4,\gamma}}}\biggr)\\
&&\qquad\quad\hspace*{38pt}{} +
|I_j^{\mathrm{side}}| \sqrt{\frac{\lambda_{4, \gamma} }{
\lambda_{2, \gamma}}} \phi\biggl(\frac{a_j C_j}{\sqrt{\lambda_{4,
\gamma}}}\biggr) +
|I_j^{\mathrm{mode}}| \sqrt{\frac{\lambda_{6, \gamma} }{
\lambda_{4, \gamma}}} \phi\biggl(\frac{a_j D_j}{\sqrt{\lambda_{4,
\gamma}}}\biggr) \Biggr] \\
&&\qquad\ge
1 - \frac{J}{L} L \biggl[
5 - 4\Phi\biggl(\frac{a C}{\sqrt{\lambda_{2,\gamma}}}\biggr)
- \Phi\biggl(\frac{a D}{\sqrt{\lambda_{4,\gamma}}}\biggr)\biggr]\\
&&\qquad\quad{}-
L \sqrt{\frac{\lambda_{4, \gamma} }{ \lambda_{2, \gamma}}} \phi
\biggl(\frac{a C}{\sqrt{\lambda_{4, \gamma}}}\biggr) -
L \sqrt{\frac{\lambda_{6, \gamma} }{ \lambda_{4, \gamma}}} \phi
\biggl(\frac{a D}{\sqrt{\lambda_{4, \gamma}}}\biggr).
\end{eqnarray*}
But this bound goes to 1 under condition (C2) as in part (1).

(3) The required probability to obtain exactly $J$ local maxima in the
set $\mathbb{S}_{1, \gamma} = \bigcup_{j=1}^{J} S_{j, \gamma}$ that
exceed $u$ is greater than the probability that exactly one local
maximum exceeds $u$ in each interval $I_j^{\mathrm{mode}}$. This probability
is bounded below by
\[
\P\Biggl[\bigcap_{j=1}^{J} \bigl(\#\{t \in\tilde{T}\cap I_j^{\mathrm{mode}}\dvtx y_\gamma(t) > u \} = 1
\cap\#\{t \in\tilde{T}\cap I_j^{\mathrm{side}} \} = 0 \bigr) \Biggr],
\]
but this goes to 1 by a similar argument as the one in part (2) of this lemma.

(4)
Since $\tilde{m}_{1, \gamma}/L = (\tilde{m}_{1, \gamma}/J)(J/L)$,
with $J/L \to A_1$, we need to show that $\tilde{m}_{1,\gamma}/J \to
1$ in probability. For any fixed $\eps>0$,
\[
0 \le\P\biggl(\biggl|\frac{\tilde{m}_{1, \gamma}}{J}-1 \biggr|
\ge\eps\biggr)
= \P( |\tilde{m}_{1, \gamma} - J| \ge J\eps
) \le
\P( \tilde{m}_{1, \gamma} \neq J) = 1-\P(\tilde
{m}_{1, \gamma} = J )
\]
since $\tilde{m}_{1, \gamma}$ and $J$ are integers. But the
probability to get exactly $J$ local maxima goes to 1 by part (2) of
this lemma.

(5)
By part (2) of this lemma, $\P[W_\gamma(u) = J]\to1$ in probability;
therefore, using the same arguments as in part (4) of this lemma, we
get $W_\gamma(u)/J \to1$. Now,
\[
\frac{W_\gamma(u)}{\tilde{m}_{1, \gamma}} =
\frac{W_\gamma(u)}{J}\frac{J}{\tilde{m}_{1, \gamma}}.
\]
But $\tilde{m}_{1, \gamma}/J \to1$ by part (3) of this lemma.
\end{pf}

\subsection{Strong control of FWER}
\label{appFWER}

%
\begin{lemma}
\label{lemmathreshold}
Let $\tilde{m}_{0,\gamma}$ be the number of local maxima in $\mathbb
{S}_{0,\gamma}$ as in Lem\-ma~\ref{lemmaPalm-S0}. Define the thresholds
$\tilde{v}_{\Bon} = F_\gamma^{-1}(\alpha/\tilde{m}_{0,
\gamma})$, random, and $v^*_{\Bon} =\break F_\gamma^{-1}
(\alpha/ \E[\tilde{m}_{0, \gamma}])$, deterministic. Then
$|\tilde{v}_{\Bon} - v^*_{\Bon}| \to0$ in probability as $L \to
\infty$.
\end{lemma}
\begin{pf}
By ergodicity, the weak law of large numbers gives that
%
\begin{equation}
\label{eqWLLN}
\biggl|\frac{\tilde{m}_{0,\gamma}}{L} - \E[\tilde{m}_{0,\gamma
}(0,1)]\biggr| \to0\vadjust{\goodbreak}
\end{equation}
in probability as $L \to\infty$, where $\E[\tilde{m}_{0,\gamma}(0,1)]
= \E[\tilde{m}_{0,\gamma}]/L$, given by
(\ref{eqexpected-local-maxima}), does not depend on $L$
[\citet{Cramer1967}]. Since $\log(\cdot)$ is continuous, the continuous
mapping theorem gives that
\[
\biggl|\log\frac{\tilde{m}_0}{L} - \log\frac{\E[\tilde
{m}_{0,\gamma}]}{L}
\biggr|
= \biggl|\log\frac{\tilde{m}_0}{\alpha} - \log\frac{\E[\tilde
{m}_{0,\gamma}]}{\alpha} \biggr| \to0,
\]
where we have used the additive property of the logarithm.

Define now the monotone increasing function $\psi_\gamma(x) =
F^{-1}_\gamma(1-e^{-x})$. The function $\psi_\gamma(x)$ is Lipschitz
continuous for all $x > 1$ because its derivative $d\psi_\gamma(x)/dx
= e^{-x}/\dot{F}_\gamma[\psi_\gamma(x)]$ is bounded for all $x >
1$. Hence, as $L \to\infty$,
\[
\biggl|\psi_\gamma\biggl(\log\frac{\tilde{m}_{0,\gamma}}{\alpha
}\biggr) -
\psi_\gamma\biggl(\log\frac{\E[\tilde{m}_{0,\gamma}]}{\alpha
}\biggr)
\biggr|
= |\tilde{v}_{\Bon} - v^*_{\Bon}| \to0.
\]
\upqed\end{pf}
\begin{pf*}{Proof of Theorem \ref{thmFWER}}
Let $\tilde{m}_{0, \gamma} \le\tilde{m}$ be the number of local
maxima in the set $\mathbb{S}_{0, \gamma}$ as in Lemma \ref
{lemmathreshold}, and let $\tilde{v}_{\Bon} = F_\gamma^{-1}
(\alpha/\tilde{m}_{0, \gamma}) \le\tilde{u}_{\Bon}$. Then
$\FWER(\tilde{u}_{\Bon}) \le\FWER(\tilde{v}_{\Bon})$. Further,
the bound $\FWER(\tilde{v}_{\Bon})$ is the probability of obtaining
at least one local maximum greater than $\tilde{v}_{\Bon}$ in
$\mathbb{S}_0 = \mathbb{S}_{0,\gamma} \cup\mathbb{T}_\gamma$,
which is less than the probability of obtaining at least one local
maximum greater than $\tilde{v}_{\Bon}$ in $\mathbb{S}_{0,\gamma}$
or at least one local maximum in~$\mathbb{T}_\gamma$.
%
\begin{equation}
\label{eqsFWER2}
\FWER(\tilde{u}_{\Bon}) \le\P[V_\gamma(\tilde{v}_{\Bon
})\ge1] +
\P(\# \{t\in\tilde{T} \cap\mathbb{T}_\gamma\}\ge1 ),
\end{equation}
where $V_\gamma(u) = \#\{t \in\tilde{T} \cap\mathbb{S}_{0,\gamma
}\dvtx y_\gamma(t)>u \}$ as in Lemma \ref{lemmaPalm-S0}.

The second probability in (\ref{eqsFWER2}) goes to zero by Lemma \ref
{lemmaunique-max}, part (1). To bound the first probability in (\ref
{eqsFWER2}), write
\[
\P[V_\gamma(\tilde{v}_{\Bon})\ge1] = \P\Bigl\{
\tilde{T} \cap\mathbb{S}_{0,\gamma} \ne\varnothing\mbox{ and }
\max_{t \in\tilde{T}} y_\gamma(t) > (\tilde{v}_{\Bon} - v^*_{\Bon
}) + v^*_{\Bon} \Bigr\},
\]
where $v^*_{\Bon} = F_\gamma^{-1}(\alpha/\E[\tilde{m}_{0,
\gamma}])$ is deterministic. For any two random variables~$X$,
$Y$ and any two constants $c$, $\eps$:
$\P(X > Y + c) \le \P(X > c - \eps) + \P(|Y| > \eps)$. Applying this
inequality with $X = \max_{t \in\tilde{T}} y_\gamma(t)$, $Y =
\tilde{v}_{\Bon} - v^*_{\Bon}$ and \mbox{$c = v^*_{\Bon}$},
%
\begin{eqnarray}
\label{eqsFWER3}
\P[V_\gamma(\tilde{v}_{\Bon})\ge1] &\le&\P
[V_\gamma(v^*_{\Bon} - \eps)\ge1]\nonumber\\[-8pt]\\[-8pt]
&&{} +
\P\{\tilde{T} \cap\mathbb{S}_{0,\gamma} \ne\varnothing
\mbox{ and } |\tilde{v}_{\Bon} - v^*_{\Bon}| > \eps\}.\nonumber
\end{eqnarray}

The second summand goes to 0 in probability as $L \to\infty$ by Lemma
\ref{lemmathreshold}. For the first summand, Lemma \ref{lemmaPalm-S0}
with level $v^*_{\Bon}-\eps$ gives that
\begin{eqnarray*}
\P[V_\gamma(v^*_{\Bon} - \eps)\ge1]
&\le&\E[V_\gamma(v^*_{\Bon} - \eps)] = \E[\tilde{m}_{0,\gamma}]
F_\gamma(v^*_{\Bon}-\eps)\\
&=& \alpha\frac{ F_\gamma(v^*_{\Bon} - \eps)}{ F_\gamma(v^*_{\Bon})},
\end{eqnarray*}
but the last fraction goes to 1 as $L \to\infty$. Replacing in (\ref
{eqsFWER3}) and (\ref{eqsFWER2}) gives the result.
\end{pf*}

\subsection{Control of FDR}
\label{appFDR}

\begin{lemma}
\label{lemmamFDR}
For any nonnegative integer random variables $V$, $W$ and fixed
positive integer $J$,
\[
\E\biggl(\frac{V}{V+W} \biggr) \le\P(W \le J-1) + \frac{\E
[V]}{\E[V]+J }.
\]
\end{lemma}
\begin{pf}
\begin{eqnarray*}
\E\biggl(\frac{V}{V+W} \biggr)
&=& \sum_{v=0}^\infty\sum_{w=0}^{J-1} \biggl(\frac{v}{v+w}
\biggr)\P(V=v, W=w ) \\[-3pt]
&&{} +
\sum_{v=0}^\infty\sum_{w=J}^\infty\biggl(\frac{v}{v+w} \biggr)\P
(V=v, W=w )
\\[-3pt]
&\le&\sum_{w=0}^{J-1}\sum_{v=0}^\infty\P(V=v, W=w )\\[-3pt]
&&{} +
\sum_{v=0}^\infty
\sum_{w=J}^\infty\biggl(\frac{v}{v+J} \biggr)\P(V=v, W=w
)
\\[-3pt]
&\le&\P( W \le J-1 ) + \E\biggl(\frac{V}{V+J} \biggr)\\[-3pt]
&\le&\P( W \le J-1 ) + \frac{\E(V)}{\E(V)+J}.
\end{eqnarray*}
The last inequality holds by Jensen's inequality, since $V/(V+J)$ is a concave
function of $V$ for $V \ge0$ and $J \ge1$.
\end{pf}
%
%
\begin{pf*}{Proof of Theorem \ref{thmFDR}}
Let $\tilde{G}(u) = \#\{t\in\tilde{T}\dvtx y_\gamma(t)>u\}/\#\{t\in
\tilde{T}\}$ be the empirical marginal right cdf of $y_\gamma(t)$
given $t \in\tilde{T}$. Then the BH threshold~$\tilde{u}_{\BH}$~%
(\ref{eqthresh-BH-random}) satisfies $\alpha\tilde{G}(\tilde
{u}_{\BH}) = k\alpha/\tilde{m} = F_\gamma(\tilde{u}_{\BH})$, so
$\tilde{u}_{\BH}$ is the largest $u$ that solves the equation
%
\begin{equation}
\label{eqFDRthreshold}
\alpha\tilde{G}(u) = F_\gamma(u).
\end{equation}
The strategy is to solve equation (\ref{eqFDRthreshold}) in the limit
when $L, a \to\infty$. We first find the limit of $\tilde{G}(u)$.
Letting $V_\gamma(u) = \#\{t \in\tilde{T} \cap\mathbb{S}_{0,\gamma
}\dvtx y_\gamma(t)>u \}$ as in Lemma \ref{lemmaPalm-S0} and $W_\gamma(u)
= \#\{t \in\tilde{T} \cap\mathbb{S}_{1,\gamma}\dvtx y_\gamma(t)>u \}
$, so that $R_\gamma(u) = V_\gamma(u)+ W_\gamma(u)$, write
%
\begin{equation}
\label{eqecdf}
\tilde{G}(u) = \frac{R_\gamma(u)}{\tilde{m}} =
\frac{V_\gamma(u)} {\tilde{m}_{0, \gamma}}
\frac{\tilde{m}_{0, \gamma}}{\tilde{m}_{0, \gamma} +\tilde
{m}_{1,\gamma}}
+ \frac{W_\gamma(u)}{\tilde{m}_{1, \gamma}}
\frac{\tilde{m}_{1, \gamma}}{\tilde{m}_{0, \gamma} +\tilde
{m}_{1,\gamma}}.
\end{equation}
By the weak law of large numbers (\ref{eqWLLN}) and Lemma \ref
{lemmaunique-max}, part (3),
\[
\frac{\tilde{m}_{0, \gamma}}{\tilde{m}_{0, \gamma} +\tilde{m}_{1,
\gamma}} = \frac{\tilde{m}_{0, \gamma}/L}{\tilde{m}_{0, \gamma}/L
+\tilde{m}_{1, \gamma}/L} \to
\frac{\E[\tilde{m}_{0, \gamma}(0,1)]}{\E[\tilde{m}_{0,\gamma
}(0,1)] + A_1}\vadjust{\goodbreak}
\]
as $L\to\infty$, where the expectation is given by (\ref
{eqexpected-local-maxima}). In addition we have the results of Lemma
\ref{lemmaPalm-S0} and Lemma \ref{lemmaunique-max}, parts (4) and
(5). Replacing these three limits in (\ref{eqecdf}), we obtain
\[
\tilde{G}(u) \to F_\gamma(u)\frac{\E[\tilde{m}_{0, \gamma
}(0,1)]}{\E[\tilde{m}_{0, \gamma}(0,1)] + A_1} + \frac{A_1}{\E
[\tilde{m}_{0, \gamma}(0,1)] + A_1}.
\]
Now replacing $\tilde{G}(u)$ by its limit in (\ref{eqFDRthreshold}),
and solving for $u$ gives the deterministic solution
%
\begin{equation}
\label{equBH^*}
F_\gamma(u^*_{\BH}) = \frac{\alpha A_1}{A_1 + \E[\tilde{m}_{0,
\gamma}(0,1)](1-\alpha)}.
\end{equation}

The FDR at the threshold $u^*_{\BH}$ is bounded by Lemma \ref
{lemmamFDR} by
%
\begin{eqnarray}
\label{eqFDR-u*}
\FDR(u^*_{\BH}) &\le& \P\bigl(W(u^*_{\BH})\le J-1\bigr) +
\frac{ \E[ V(u^*_{\BH})]}{ \E[V(u^*_{\BH
})] + J } \nonumber\\
&=&\P\bigl(W(u^*_{\BH})\le J-1\bigr)\\
&&{} + \frac{\E[ V_\gamma(u^*_{\BH
})]
+ \E[ \#\{t \in\tilde{T} \cap\mathbb{T}_\gamma\dvtx
y_\gamma(t)>u^*_{\BH} \} ] }{ \E[V_\gamma(u^*_{\BH
})]
+ \E[ \#\{t \in\tilde{T} \cap\mathbb{T}_\gamma\dvtx
y_\gamma(t)>u^*_{\BH} \} ] + J },\nonumber
\end{eqnarray}
where we have split $V_\gamma(u^*_{\BH})$ into the reduced null
region $\mathbb{S}_{0, \gamma}$ and the transition region $\mathbb
{T}_\gamma= \mathbb{S}_0 \setminus\mathbb{S}_{0, \gamma}$. Under
condition (C2), Lemma \ref{lemmaunique-max}, part (1), gives
%
\begin{equation}
\label{eqlim2}
0 \le\E[ \#\{t\in\tilde{T}\cap\mathbb{T}_\gamma\dvtx y_\gamma
(t)>u^*_{\BH} \} ]
\le\E[ \#\{t\in\tilde{T}\cap\mathbb{T}_\gamma\} ]
\to0.
\end{equation}
By Lemma \ref{lemmaPalm-S0}, the remaining terms of the last fraction
in (\ref{eqFDR-u*}) can be written as
\begin{eqnarray*}
\frac{ \E[ V_\gamma(u^*_{\BH})]}{ \E[
V_\gamma(u^*_{\BH})]+J} &=&
\frac{F_\gamma(u^*_{\BH})\E[\tilde{m}_{0, \gamma}(0,1)]L}
{F_\gamma(u^*_{\BH})\E[\tilde{m}_{0, \gamma}(0,1)]L + J} \\
&=&
\frac{F_\gamma(u^*_{\BH})\E[\tilde{m}_{0, \gamma}(0,1)]}
{F_\gamma(u^*_{\BH})\E[\tilde{m}_{0, \gamma}(0,1)] + J/L}.
\end{eqnarray*}
Since $u^*_{\BH}$ solves (\ref{equBH^*}), for $L \to\infty$ such
that $J/L \to A_1$, the above expression tends to
%
\begin{equation}
\label{eqlim1}
\frac{\alpha\E[\tilde{m}_{0, \gamma}(0,1)]}{\alpha\E[\tilde
{m}_{0, \gamma}(0,1)] + A_1 + (1-\alpha)\E[\tilde{m}_{0, \gamma}(0,1)]}
= \alpha\frac{\E[\tilde{m}_{0, \gamma}(0,1)]}{\E[\tilde{m}_{0,
\gamma}(0,1)] + A_1} \le\alpha.\hspace*{-28pt}
\end{equation}
Combining equations (\ref{eqlim2}), (\ref{eqlim1}) and Lemma \ref
{lemmaunique-max}, part (3), in (\ref{eqFDR-u*}), we obtain $\lim\sup
\FDR(u^*_{\BH}) \le\alpha$.

Recall that the BH threshold $\tilde{u}_{\BH}$ solves equation (\ref
{eqFDRthreshold}), and $u^*_{\BH}$ satisfies~(\ref{equBH^*}), where
the empirical marginal distribution, $\tilde{G}(u)$, is replaced by
its limit. Since $F_\gamma(t)$ is continuous, $F_\gamma(\tilde
{u}_{\BH})\to F_\gamma(u^*_{\BH})$, leading to\break \mbox{$\lim\sup\FDR
(\tilde{u}_{\BH}) \le\alpha$}.
\end{pf*}

\subsection{Power}\vspace*{-3pt}
\label{apppower}

\begin{lemma}
\label{lemmafwer-fdr-threshold}
For any $j = 1,\ldots,J$, let $t$ be any interior point of the support
$S_j$ of peak $j$. Under conditions \textup{(C1)} and \textup{(C2)},
\[
u^*_{\Bon}/[a_j h_{j,\gamma}(t)] \to0, \qquad
u^*_{\BH}/[a_j h_{j,\gamma}(t)] \to0
\]
in probability, where $u^*_{\Bon}$ and $u^*_{\BH}$ are given by
(\ref{eqthresh-Bon-fixed}) and (\ref{eqthresh-BH-fixed}), respectively.\vspace*{-3pt}
\end{lemma}
\begin{pf}
(1) From (\ref{eqdistr}), for $u > \sigma_\gamma$, $F_\gamma(u)$ is
bounded above and below by
%
\begin{equation}
\label{eqF-bounds}
\frac{C_1}{2} \phi\biggl(\frac{u}{\sigma_\gamma}\biggr)
< F_\gamma(u) <
(C_1 + 1) \phi\biggl(\frac{u}{\sigma_\gamma}\biggr),\qquad
C_1 = \sqrt{\frac{2\pi\lambda_{2,\gamma}^2}{\lambda_{4,\gamma}
\sigma^2_\gamma}},
\end{equation}
where the lower bound was obtained using $\Phi(x) > 1/2$ for $x>1$,
and the upper bound used the fact that $\sqrt{\lambda_{4,\gamma
}/\Delta} \ge1/\sigma_\gamma$ and $1-\Phi(x) < \phi(x)/x$ for
$x>1$. Let $v = F_\gamma(u)$. Inverting the bounds in (\ref
{eqF-bounds}) we obtain
%
\begin{equation}
\label{eqinverse-bounds}
2\sigma^2_\gamma\biggl(\log{\frac{C_1}{2\sqrt{2\pi}}} - \log
v\biggr)
< u^2 <
2\sigma^2_\gamma\biggl(\log{\frac{C_1+1}{\sqrt{2\pi}}} - \log
v\biggr).
\end{equation}
Applying these inequalities to $v^* = F_\gamma(u^*_{\Bon})$ and $w =
F_\gamma[a_j h_{j,\gamma}(t)]$ gives that
\[
0 \le\frac{(u^*_{\Bon})^2}{[a_j h_{j,\gamma}(t)]^2} <
\frac{\log[(C_1+1)/\sqrt{2\pi}] - \log(v^*)}{\log[C_1/(2\sqrt
{2\pi})] - \log(w)}.
\]
Applying L'H\^{o}pital's rule, the limit of the above fraction when
$v^*$ and $w$ go to zero is the same as the limit of $w/v^*$. But this
limit is zero because, by the upper bound in (\ref{eqF-bounds}) and
(\ref{eqthresh-Bon-fixed}),
\[
\frac{F_\gamma[a_j h_{j,\gamma}(t)]}{F_\gamma(u^*_{\Bon})} <
(C_1+1)\frac{A_1 + \E[\tilde{m}_{0,\gamma}(0,1)]}{\alpha} L\phi
\biggl(\frac{a_j h_{j,\gamma}(t)}{\sigma_\gamma}\biggr),
\]
which goes to zero by the lemma's conditions.

(2) The FDR threshold $u^*_{\BH}$ (\ref{eqthresh-BH-fixed}) is
bounded, so the result is immediate.\vspace*{-3pt}
\end{pf}
%
%
\begin{pf*}{Proof of Theorem \ref{thmpower}}
For any threshold $u$, the detection power\break $\Power(u)$~(\ref
{eqpower}) is greater than $\E[W_\gamma(u)] / J \ge\P[W_\gamma(u) =
J]$. But this probability goes to~1 by Lemma \ref{lemmaunique-max},
part (3), particularly for the deterministic thresholds~$u^*_{\Bon}$
and $u^*_{\BH}$. It was shown in the proofs of Theorems \ref{thmFWER}
and \ref{thmFDR} that the gap between the deterministic thresholds and
the random thresholds~$\tilde{u}_{\Bon}$ and~$\tilde{u}_{\BH}$
narrows to zero asymptotically. Therefore the power for these
thresholds goes to 1 as well.\vspace*{-3pt}
\end{pf*}

\section*{Acknowledgments}

The authors thank Pablo Jadzinsky for providing the neural recordings
data, as well as Igor Wigman, Felix Abramovich and Yoav Benjamini for
helpful discussions. The authors also thank the Editor, Associate
Editor and referees for their handling of the manuscript and their
useful suggestions.\vadjust{\goodbreak}


%

\printaddresses

\end{document}